\documentclass{amsart}    

\usepackage[francais]{babel}
\usepackage{amsmath}
\usepackage[latin1]{inputenc}
\usepackage{amsfonts}
\usepackage{amssymb}
\usepackage{graphicx,epsfig}
\usepackage{amsthm}
\usepackage{amscd}
\usepackage[all]{xy}
\usepackage[T1]{fontenc}
\usepackage[francais]{babel}

\usepackage[latin1]{inputenc}

\usepackage[T1]{fontenc}

\usepackage{amsmath}

\usepackage{amsfonts}

\usepackage{amssymb}

\usepackage{amsthm}

\usepackage{graphics}

\newcommand{\mk}{\medskip}

\newcommand{\ZZ}{\mathbb{Z}}

\newcommand{\CC}{\mathbb{C}}

\newcommand{\NN}{\mathbb{N}}

\newcommand{\Glie}{\mathfrak{g}}

\newcommand{\Yim}{\mathcal{Y}}

\newcommand{\demo}{\noindent {\it \small D\'emonstration:}\quad}

\renewcommand{\NN}{\ensuremath{\mathbb{N}}}

\renewcommand{\CC}{\ensuremath{\mathbb{C}}}

\newcommand{\U}{\mathcal{U}}

\newtheorem{defi}{D\'{e}finition}

\newtheorem{thm}{Th\'{e}or\`{e}me}

\newtheorem{cor}{Corollaire}

\newtheorem{prop}{Proposition}

\newtheorem{lem}{Lemme}

\title{$t$-analogues des op\'erateurs d'\'ecrantage associ\'es aux $q$-caract\`eres}
\author{David Hernandez}
\address{David Hernandez: \'Ecole Normale Sup\'erieure - DMA, 45, Rue d'Ulm F-75230 PARIS,
  Cedex 05  FRANCE}
\email{David.Hernandez@ens.fr\\ URL: http://www.dma.ens.fr/$\sim$dhernand}
\begin{document}

\begin{abstract} Nous proposons des op\'erateurs d'\'ecrantage pour la th\'eorie des $q,t$-caract\`eres de Nakajima (\cite{Naa}, \cite{Nab}), analogues aux op\'erateurs d'\'ecrantage de Frenkel et Reshetikhin relatifs \`a leur th\'eorie des $q$-caract\`eres pour les repr\'esentations de dimension finie des alg\`ebres affines quantifi\'ees \cite{Fre}, avec en particulier les m\^emes propri\'et\'es de sym\'etrie. La th\'eorie de Nakajima utilise des anneaux non-commutatifs et nous aurons ainsi \`a consid\'erer des bimodules adapt\'es \`a ces structures. Notre construction \'etant purement alg\'ebrique et s'appuyant sur la d\'efinition combinatoire des $q,t$-caract\`eres, elle est \'etendue au cas non-simplement laç\'e.
\\
\\
\begin{center}$t$-ANALOGS OF SCREENING OPERATORS RELATED TO $q$-CHARACTERS\end{center}
ABSTRACT. Frenkel and Reshetikhin introduced screening operators related to $q$-characters of finite dimensional representations of quantum affine algebras \cite{Fre}. We propose $t$-analogs of screening operators related to Nakajima's $q,t$-characters (\cite{Naa}, \cite{Nab}) with the same properties of symmetry. Nakajima introduced non-commutative rings, so we propose bimodules. Our construction uses only combinatorial definition of $q,t$-characters, and therefore can be extended to the non-simply laced case.
\end{abstract}

\maketitle

For convenience of the reader we give an english translation of the introduction:

\section*{Introduction}

\noindent Let $q\in\CC^*$ such that $q$ is not a root of unity.

\noindent In the case of semi-simple Lie algebras $\Glie$, the structure of the Grothendieck ring $\text{Rep}(\U_q(\Glie))$ of finite dimensional representations of the quantum algebra $\U_q(\Glie)$ is well understood, see \cite{Ro}. It is analogous to the classic case $q=1$. In particular we have ring isomorphisms:
$$\text{Rep}(\U_q(\Glie))\simeq \text{Rep}(\Glie)\simeq \ZZ [\Lambda]^W\simeq \ZZ[T_1,...,T_n]$$ 
deduced from the injective homomorphism of characters $\chi$:
$$\chi(V)=\underset{\lambda\in\Lambda}{\sum}\text{dim}(V_{\lambda})\lambda$$
where $V_{\lambda}$ are weight spaces of a representation $V$ and $\Lambda$ is the set of weight of $V$.

\noindent For the general case of Kac-Moody algebras the picture is less clear. In the affine case $\U_q(\hat{\Glie})$, Frenkel and Reshetikhin \cite{Fre}, motivated by the theory of deformed $W$-algebras, have recently introduced an injective ring homomorphism of $q$-characters:
$$\chi_q:\text{Rep}(\U_q(\hat{\Glie}))\rightarrow \ZZ[Y_{i,a}^{\pm 1}]_{1\leq i\leq n,a\in\CC^*}=\Yim$$
The construction of $\chi_q$ uses the universal $R$-matrix. The homomorphism $\chi_q$ allows to understand the ring $\text{Rep}(\U_q(\hat{\Glie}))\simeq\ZZ[X_{i,a}]_{i\in I,a\in\CC^*}$. The classical limit $q\rightarrow 1$ of $\chi_q$ is the usual homomorphism of characters. In fact $\chi_q$ gives informations about the decomposition in Jordan subspaces for a class $(\phi^{\pm}_{i,m})_{m\in \ZZ,i\in I}$ of commutative elements of $\U_q(\hat{\Glie})$:
$$\chi_q(V)=\underset{\gamma}{\sum}\text{dim}(V_{(\gamma)})\underset{i\in I}{\prod}\underset{r=1..k_{\gamma i}}{\prod}Y_{i,a_{\gamma ir}}\underset{s=1..l_{\gamma i}}{\prod}Y_{i,b_{\gamma ir}}^{-1}$$
where $V_{(\gamma)}$ is the Jordan subspace of weight $\gamma$:
$$\underset{m\geq 0}{\sum}\gamma_{i,\pm m}^{\pm}u^{\pm m}=\gamma_i^{\pm}(u)=q_i^{k_{\gamma i}-l_{\gamma i}}\frac{Q_i(uq_i^{-1})R_i(uq_i)}{Q_i(uq_i)R_i(uq_i^{-1})}$$
where $a_{\gamma ir}$ are roots of the polynomial $Q_i$ and $b_{\gamma ir}$ roots of the polynomial $R_i$.
\noindent The homomorphism of $q$-characters has a symmetry property analogous to the classic action of the Weyl group $\text{Im}(\chi)=\ZZ[\Lambda]^W$: Frenkel and Reshetikhin defined $n$ screnning operators (with $A_{i,a}\in \Yim$ monomials):
$$S_i:\ZZ[Y_{i,a}^{\pm 1}]_{a\in\CC^*,i\in I}=\Yim\rightarrow\underset{a\in\CC^*}{\bigoplus}\Yim.S_{i,a} / \underset{a\in\CC^*}
{\sum}\Yim.(S_{i,aq_{i}^2}-A_{i,aq_i}.S_{i,a})$$
There is a leibnitz rule ($S_i(UV)=US_i(V)+VS_i(U)$), and $S_i(Y_a)=Y_a.S_a$. They conjectured:
$$\text{Im}(\chi_q)=\underset{i\in I}{\bigcap}\text{Ker}(S_i)$$
They proved it in the case $sl_2$ \cite{Fre}, and Frenkel, Mukhin proved it in the general case \cite{Fre2}.

\noindent These operators give informations about the combinatorial structure of $q$-characters. For example:
$$\text{Rep}(\U_q(\hat{\Glie}))\simeq \text{Im}(\chi_q)=\mathfrak{K}=\underset{i\in I}{\bigcap}(\ZZ[Y^{\pm}_{j,a}]_{j\neq i,a\in\CC^*}\ZZ[Y_{i,b}+Y_{i,b}A^{-1}_{i,bq_i}]_{b\in\CC^*})$$

\noindent In the $ADE$ case Nakajima introduced $t$-analogs of $q$-characters \cite{Naa}, \cite{Nab}. He defined maps $\chi_{q,t}$ et $\hat{\chi}_{q,t}$ from $\text{Rep}(\U_q(\hat{\Glie}))\otimes_{\ZZ}\ZZ[t,t^{-1}]$ to polynomial rings respectively $\Yim_t=\ZZ[Y_{i,a}^{\pm},t^{\pm}]_{i\in I,a\in\CC^*}$ and 
\\$\hat{\Yim}_t=\ZZ[V_{i,a},W_{i,a},t^{\pm}]_{i\in I,a\in\CC^*}$ . From representation theory point of view, it gives more informations about Jordan subspaces. He introduces a new non-commutative multiplication $*$ on $\hat{\Yim}_t$ and gives a combinatorial axiomatic definition of $q,t$-characters. The existence is non-trivial and is proved with geometric approach. 

\noindent In this paper we propose $t$-analogs of screening operators $\hat{S}_{i,t}, S_{i,t}$ related to applications $\chi_{q,t}$ and $\hat{\chi}_{q,t}$.

\noindent This article is organized as follows. In section 2 we recall the fundamental symmetry property of Frenkel, Reshetikhin's screnning operators and some results on Nakajima's ring $\hat{\Yim}_t$. In section 3 we define operators $\hat{S}_{t,i}^l$. We introduce a bimodule structure such that we have a Leibnitz rule. In section 4 we define $t$-analogs of screening operators $\hat{S}_{t,i}$. In the $ADE$ case we give an interpretation related to Nakajima's multiplication $*$. These operators verify the expected symmetry property in theorem \ref{justi}: 
$$\underset{i\in I}{\bigcap}\text{Ker} (\hat{S}_i)=\hat{\mathfrak{K}}\supseteq \text{Im}(\hat{\chi}_{q,t})$$
In section 5 we define operators $S_{i,t}$ related to the ring $\Yim_t$. The diagram:
$$\begin{array}{rcccl}
\hat{\Yim}_t&\stackrel{\hat{S}_{i,t}}{\longrightarrow}&\hat{\Yim}_{t,i}\\
\hat{\Pi}_t\downarrow &&\downarrow&\hat{\Pi}_{t,i}\\
\Yim_t&\stackrel{S_{t,i}}{\longrightarrow}&\Yim_{t,i}\\
\Pi_t\downarrow &&\downarrow&\Pi_{t,i}\\
\Yim&\stackrel{S}{\longrightarrow}&\Yim_i\\
\end{array}$$ 
is commutative. We have a symmetry property in theorem \ref{ksimple}:
$$\underset{i\in I}{\bigcap}\text{Ker} (S_i)=\mathfrak{K}=\text{Im}(\chi_{q,t})$$
where $\chi_{q,t}$ is $\tilde{\chi}_{q,t}$ in \cite{Naa}. In section 6 we construct involutions analog to the Nakajima's one. 

\noindent The construction uses a bimodule structure on the free left module $\underset{a\in\CC^*}{\bigoplus}\hat{\Yim}_tS_{i,a}$, and the $t$-analog of 
\\$\underset{a\in\CC^*}{\sum}\Yim.(S_{i,aq_{i}^2}-A_{i,aq_i}.S_{i,a})$ is a subbimodule. The bimodule structure is given for any $m\in\hat{\Yim}_t$ monomial by:
$$S_{i,a}m=t^{2u_{i,a}(m)}mS_{i,a}$$

\section{Introduction}

\noindent Dans ce qui suit $q\in\CC^*$ est suppos\'e ne pas \^etre une racine de l'unit\'ee.

\noindent Dans le cas d'une alg\`ebre de Lie semi-simple $\Glie$, la structure de l'anneau de Grothendieck $\text{Rep}(\U_q(\Glie))$ des repr\'esentations de dimensions finie de l'alg\`ebre semi-simple quantifi\'ee $\U_q(\Glie)$ est bien comprise, voir \cite{Ro}. En fait on a pu montrer qu'elle est tout \`a fait analogue \`a celle du cas classique $q=1$ d\'ej\`a bien connu. On a en particulier des isomorphismes d'anneaux:
$$\text{Rep}(\U_q(\Glie))\simeq \text{Rep}(\Glie)\simeq \ZZ [\Lambda]^W\simeq \ZZ[T_1,...,T_n]$$ 
construits \`a partir d'un morphisme de caract\`ere $\chi$ tel que pour une repr\'esentation $V$ de sous-espaces de poids $V_{\lambda}$:
$$\chi(V)=\underset{\lambda\in\Lambda}{\sum}\text{dim}(V_{\lambda})\lambda$$
o\`u $\Lambda$ d\'esigne l'ensemble des poids de $V$.

\noindent Par contre la quantification  modifie la th\'eorie des repr\'esentations lorsqu'on s'int\'eresse au cas g\'en\'eral des alg\`ebres de Kac-Moody. Dans le cas affine $\U_q(\hat{\Glie})$, Frenkel et Reshetikhin \cite{Fre}, motiv\'es par la th\'eorie des $W$-alg\`ebres d\'eform\'ees, ont r\'ecemment introduit un morphisme d'anneau injectif, dit de $q$-caract\`eres, \`a valeurs dans un anneau de polyn\^omes de Laurent:
$$\chi_q:\text{Rep}(\U_q(\hat{\Glie}))\rightarrow \ZZ[Y_{i,a}^{\pm 1}]_{1\leq i\leq n,a\in\CC^*}=\Yim$$
La construction de $\chi_q$ repose sur l'existence d'une $R$-matrice universelle. L'application $\chi_q$ permet de comprendre l'anneau $\text{Rep}(\U_q(\hat{\Glie}))\simeq\ZZ[X_{i,a}]_{i\in I,a\in\CC^*}$ et lorsqu'on regarde la limite classique $q=1$ on retrouve l'application de caract\`eres usuelle. En fait cette application prend en compte la d\'ecomposition en sous-espaces de Jordan pour une certaine famille commutante $(\phi^{\pm}_{i,m})_{m\in \ZZ,i\in I}$ d'\'el\'ements de $\U_q(\hat{\Glie})$:
$$\chi_q(V)=\underset{\gamma}{\sum}\text{dim}(V_{(\gamma)})\underset{i\in I}{\prod}\underset{r=1..k_{\gamma i}}{\prod}Y_{i,a_{\gamma ir}}\underset{s=1..l_{\gamma i}}{\prod}Y_{i,b_{\gamma ir}}^{-1}$$
o\`u $V_{(\gamma)}$ d\'esigne le sous-espace de Jordan de $V$ de poids:
$$\underset{m\geq 0}{\sum}\gamma_{i,\pm m}^{\pm}u^{\pm m}=\gamma_i^{\pm}(u)=q_i^{k_{\gamma i}-l_{\gamma i}}\frac{Q_i(uq_i^{-1})R_i(uq_i)}{Q_i(uq_i)R_i(uq_i^{-1})}$$
avec $a_{\gamma ir}$ les racines du polyn\^ome $Q_i$ et $b_{\gamma ir}$ les racines du polyn\^ome $R_i$.

\noindent Le morphisme de $q$-caract\`ere v\'erifie une propri\'et\'e de sym\'etrie analogue au cas classique de l'action du groupe de Weyl qui veut $\text{Im}(\chi)=\ZZ[\Lambda]^W$. En effet Frenkel et Reshetikhin ont d\'efini $n$ op\'erateurs dits d'\'ecrantage (avec $A_{i,a}\in \Yim$ mon\^omes):
$$S_i:\ZZ[Y_{i,a}^{\pm 1}]_{a\in\CC^*,i\in I}=\Yim\rightarrow\underset{a\in\CC^*}{\bigoplus}\Yim.S_{i,a} / \underset{a\in\CC^*}
{\sum}\Yim.(S_{i,aq_{i}^2}-A_{i,aq_i}.S_{i,a})$$
qui sont des d\'erivations ($S_i(UV)=US_i(V)+VS_i(U)$), v\'erifiant $S_i(Y_a)=Y_a.S_a$, et ont conjectur\'e:
$$\text{Im}(\chi_q)=\underset{i\in I}{\bigcap}\text{Ker}(S_i)$$
Ils l'ont montr\'e dans le cas $sl_2$ \cite{Fre} puis Frenkel et Mukhin ont obtenu le r\'esultat dans le cas g\'en\'eral \cite{Fre2}.

\noindent Ces op\'erateurs permettent de comprendre la structure combinatoire des $q$-caract\`eres. Par exemple:
$$\text{Rep}(\U_q(\hat{\Glie}))\simeq \text{Im}(\chi_q)=\mathfrak{K}=\underset{i\in I}{\bigcap}(\ZZ[Y^{\pm}_{j,a}]_{j\neq i,a\in\CC^*}\ZZ[Y_{i,b}+Y_{i,b}A^{-1}_{i,bq_i}]_{b\in\CC^*})$$

\noindent Dans le cas o\`u $\Glie$ est de type $ADE$, Nakajima a raffin\'e la th\'eorie en introduisant un $t$-analogue des $q$-caract\`eres \cite{Naa}, \cite{Nab}. Il consid\`ere des applications $\chi_{q,t}$ et $\hat{\chi}_{q,t}$ de $\text{Rep}(\U_q(\hat{\Glie}))\otimes_{\ZZ}\ZZ[t,t^{-1}]$ vers les anneaux de polyn\^omes respectivement $\Yim_t=\ZZ[Y_{i,a}^{\pm},t^{\pm}]_{i\in I,a\in\CC^*}$ et $\hat{\Yim}_t=\ZZ[V_{i,a},W_{i,a},t^{\pm}]_{i\in I,a\in\CC^*}$. D'un point de vue des repr\'esentations, elles permettent de mieux comprendre la structure de chaque sous-espace de Jordan. Au passage il introduit une nouvelle multiplication $*$ sur $\hat{\Yim}_t$ qui n'est pas commutative. Il donne une d\'efinition axiomatique combinatoire des $q,t$-caract\`eres. L'existence est non-triviale et est prouv\'ee gr\^ace \`a une approche g\'eom\'etrique.

\noindent Nous proposons dans cet article des $t$-analogues des op\'erateurs d'\'ecrantage, adapt\'es aux applications $\chi_{q,t}$ et $\hat{\chi}_{q,t}$. 

\noindent Dans la deuxi\`eme partie, on rappelle la propri\'et\'e fondamentale de sym\'etrie des op\'erateurs d'\'ecrantage de Frenkel et Reshetikhin (Th\'eor\`eme \ref{frene}) ainsi que quelques r\'esultats \'el\'ementaires sur l'anneau $\hat{\Yim}_t$ de Nakajima. On d\'efinit dans la troisi\`eme partie les op\'erateurs $\hat{S}^l_{t,i}$ qui peuvent \^etre interpr\'et\'es comme des d\'erivations pour la multiplication usuelle et une certaine structure de bimodule. Dans la quatri\`eme partie on d\'efinit les $t$-analogues des op\'erateurs d'\'ecrantage $\hat{S}_{t,i}$ qui dans le cas o\`u $\Glie$ est de type $ADE$ peuvent \^etre interpr\'et\'es comme des d\'erivations en utilisant la loi $*$ de Nakajima. Ces op\'erateurs v\'erifient la propri\'et\'e attendue dans le th\'eor\`eme \ref{justi}:
$$\underset{i\in I}{\bigcap}\text{Ker} (\hat{S}_i)=\hat{\mathfrak{K}}\supseteq \text{Im}(\hat{\chi}_{q,t})$$
On d\'efinit dans la cinqui\`eme partie des op\'erateurs $S_{i,t}$ pour l'anneau $\Yim_t$ rendant le diagramme commutatif: 
$$\begin{array}{rcccl}
\hat{\Yim}_t&\stackrel{\hat{S}_{i,t}}{\longrightarrow}&\hat{\Yim}_{t,i}\\
\hat{\Pi}_t\downarrow &&\downarrow&\hat{\Pi}_{t,i}\\
\Yim_t&\stackrel{S_{t,i}}{\longrightarrow}&\Yim_{t,i}\\
\Pi_t\downarrow &&\downarrow&\Pi_{t,i}\\
\Yim&\stackrel{S}{\longrightarrow}&\Yim_1\\
\end{array}$$  
avec une propri\'et\'e de sym\'etrie dans le th\'eor\`eme \ref{ksimple}:
$$\underset{i\in I}{\bigcap}\text{Ker} (S_i)=\mathfrak{K}=\text{Im}(\chi_{q,t})$$
o\`u $\chi_{q,t}$ est \'egal au $\tilde{\chi}_{q,t}$ de \cite{Naa}. Dans la sixi\`eme partie on donne la construction d'involutions analogues \`a celle de Nakajima.

\noindent Notons que la construction repose sur l'existence d'une structure de bimodule sur le module libre \`a gauche $\underset{a\in\CC^*}{\bigoplus}\hat{\Yim}_tS_{i,a}$ telle que le $t$-analogue de $\underset{a\in\CC^*}{\sum}\Yim.(S_{i,aq_{i}^2}-A_{i,aq_i}.S_{i,a})$ soit un sous-bimodule. Cette structure est caract\'eris\'ee par les relations suivantes, o\`u $m\in\hat{\Yim}_t$ est un mon\^ome:
$$S_{i,a}m=t^{2u_{i,a}(m)}mS_{i,a}$$

\section{Rappels}

Soit $\Glie$ une alg\`ebre de Lie simple. On note $n$ le rang de $\Glie$, $(C_{i,j})_{1\leq i,j\leq n}$ sa matrice de Cartan et $I=\{1,...,n\}$. On note $q_i=q^{r_i}$ comme dans \cite{Fre}.

\subsection{Op\'erateurs d'\'ecrantage \cite{Fre}, \cite{Fre2}}

On consid\`ere l'anneau:
$$\Yim=\ZZ[Y_{i,a}^{\pm 1}]_{i\in I,a\in\CC^*}$$
les $\Yim$-modules libres ($i\in I$):
$$\Yim^l_{i}=\underset{a\in\CC^*}{\bigoplus}\Yim.S_{i,a}$$
et les $\Yim$-modules $\Yim_{i}$ d\'efinis respectivement comme $\Yim$-module quotient de $\Yim^l_{i}$: 
$$\Yim_{i}=\underset{a\in\CC^*}{\bigoplus}\Yim.S_{i,a} / \underset{a\in\CC^*}
{\sum}\Yim.(S_{i,aq_i^2}-A_{i,aq_i}.S_{i,a})$$
par le sous-module $F_i=\underset{a\in\CC^*}
{\sum}\Yim.(S_{i,aq_i^2}-A_{i,aq_i}.S_{i,a})$, avec:
$$A_{i,a}=Y_{i,aq_i^{-1}}Y_{i,aq_i}\underset{j/C_{j,i}=-1}{\prod}Y_{j,a}^{-1}\underset{j/C_{j,i}=-2}{\prod}Y_{j,aq}^{-1}Y_{j,aq^{-1}}^{-1}\underset{j/C_{j,i}=-3}{\prod}Y_{j,aq^2}^{-1}Y_{j,a}^{-1}Y_{j,aq^{-2}}^{-1}$$
On a alors les op\'erateurs d'\'ecrantage:
$$S_i:\Yim\rightarrow \Yim_{i}$$
qui sont des d\'erivations pour le produit de $\Yim$:
$$S_i(U.V)=U.S_i(V)+V.S_i(U)$$
et qui v\'erifient pour $a\in\CC^*$:
$$S_i(Y_{j,a}^{\pm})=\pm \delta_{i,j}Y_{i,a}^{\pm}S_{i,a}$$
On peut d\'efinir de mani\`ere analogue $S_i^l:\Yim\rightarrow \Yim_{i}^l$.
\begin{thm}\label{frene} Le noyau de $S_i$ est le sous-anneau de $\Yim$:
$$\text{Ker}(S_i)=\mathfrak{K}_i=\ZZ[Y_{j,a}^{\pm}]_{j\neq i,a\in\CC^*}\ZZ[Y_{i,b}(1+A_{i,bq_i}^{-1})]_{b\in\CC^*}$$
\end{thm}

\subsection{L'anneau $\hat{\Yim}_t$ \cite{Nab}}

On consid\`ere \`a pr\'esent l'anneau:
$${\hat{\Yim}}_t=\ZZ[t,t^{-1},V_{i,a},W_{i,a}]_{i\in I,a\in\CC^*}$$
C'est un $\ZZ[t,t^{-1}]$-module libre de base l'ensemble des $\underset{i\in I,a\in\CC^*}{\prod}V_{i,a}^{v_{i,a}(m)}W_{i,a}^{w_{i,a}(m)}$ qu'on appelera mon\^omes.

\noindent On d\'efinit un morphisme d'anneaux:
$$\tilde{\Pi}_t:\hat{\Yim}_t\rightarrow\Yim$$
$$m=\underset{i\in I,a\in\CC^*}{\prod}V_{i,a}^{v_{i,a}(m)}W_{i,a}^{w_{i,a}(m)}\mapsto \underset{i\in I,a\in\CC^*}{\prod}Y_{i,a}^{u_{i,a}(m)}\text{ et } t\mapsto 1$$
avec pour un tel mon\^ome $m$:
$$u_{i,a}(m)=w_{i,a}(m)-v_{i,aq_i^{-1}}(m)-v_{i,aq_i}(m)$$
$$+\underset{j/C_{i,j}=-1}{\sum}v_{j,a}(m)+\underset{j/C_{i,j}=-2}{\sum}(v_{j,aq}(m)+v_{j,aq^{-1}}(m))+\underset{j/C_{i,j}=-3}{\sum}(v_{j,aq^2}(m)+v_{j,a}(m)+v_{j,aq^{-2}}(m))$$
Remarquer que l'application $\tilde{\Pi}_t$ est l'unique morphisme d'anneaux tel que:
$$\tilde{\Pi}_t(W_{i,a})=Y_{i,a}\text{ , }\tilde{\Pi}_t(V_{i,a})=A_{i,a}^{-1}\text{ , }\tilde{\Pi}_t(t)=1$$ 
On peut d\'efinir pour un mon\^ome $m=\underset{i\in I,a\in\CC^*}{\prod}Y_{i,a}^{u_{i,a}(m)}\in \Yim$ les $u_{i,a}(m)$ de mani\`ere \'evidente, et alors ces quantit\'es sont conserv\'ees par $\tilde{\Pi}_t$.

\noindent Pour $m\in\hat{\Yim}$ mon\^ome $i$-dominant, c'est \`a dire v\'erifiant $\forall a\in\CC^*, u_{i,a}(m)\geq 0$, on pose:
$$E_i(m)=m\underset{a\in\CC^*}{\prod}\underset{r_a=0..u_{i,a}(m)}{\sum}t^{r_a(u_{i,a}(m)-r_a)}\begin{bmatrix}u_{i,a}(m)\\r_a\end{bmatrix}_tV_{i,aq_i}^{r_a}$$
et on note $\hat{\mathfrak{K}}_{t,i}$ le sous $\ZZ[t,t^{-1}]$-module de $\hat{\Yim}_t$ engendr\'e par ces $E_i(m)$. On pose alors: 
$$\hat{\mathfrak{K}}_t=\underset{i\in I}{\bigcap}\hat{\mathfrak{K}}_{t,i}$$
On note $\hat{A}$ l'ensemble des mon\^omes de $\hat{\Yim}_t$, $\hat{B}_i\subset \hat{A}$ l'ensemble des mon\^omes $i$-dominants de $\hat{\Yim}_t$.

\begin{lem}\label{sdirecte}  Pour chaque $i\in I$, on a une d\'ecomposition en somme directe de $\ZZ[t,t^{-1}]$-modules:
$$\hat{\Yim}_t=\hat{\mathfrak{K}}_{t,i}\oplus \underset{m\in \hat{A}-\hat{B}_i}{\bigoplus} \ZZ[t,t^{-1}]m=(\underset{m\in \hat{B}_i}{\bigoplus}\ZZ[t,t^{-1}]E_i(m))\oplus (\underset{m\in \hat{A}-\hat{B}_i}{\bigoplus} \ZZ[t,t^{-1}]m)$$
\end{lem}

\demo

Notons d'abord que pour $m\in \hat{B}_i$, on peut \'ecrire $E_i(m)=m+f(m)$ avec 
$$f(m)=m((\underset{a\in\CC^*}{\prod}\underset{r_a=0..u_{i,a}}{\sum}t^{r_a(u_{i,a}-r_a)}\begin{bmatrix}u_{i,a}\\r_a\end{bmatrix}_tV_{i,aq_i}^{r_a})-1)$$ 
qui ne fait intervenir que des mon\^omes de $i$-poids $wt_i(m')=\underset{a\in\CC^*}{\sum}u_{i,a}(m')$ strictement inf\'erieur \`a celui de $m$.

\noindent Consid\'erons une combinaison lin\'eaire qui s'annule:
$$\underset{m\in \hat{B}_i}{\sum}\lambda_m(t)E_i(m)+\underset{m\in \hat{A}-\hat{B}_i}{\sum}\mu_m(t)m=0$$
avec les $\lambda_m(t),\mu_m(t)\in\ZZ[t,t^{-1}]$. Si on suppose qu'un des $\lambda(t)\neq 0$, soit $m_1\in \hat{B}_i$ un mon\^ome dominant de $i$-poids maximal parmi ceux qui v\'erifient $\lambda_m(t)\neq 0$. Alors le mon\^ome $m_1$ ne peut appara\^itre que dans $E_i(m_1)$ puisque si il apparaissait dans $E_i(m_2)$, le $i$-poids de $m_2$ serait strictement plus grand que le sien. Donc $\lambda_{m_1}(t)=0$, contradiction. Donc tous les $\lambda_m(t)$ sont nuls, et alors $\underset{m\in \hat{A}-\hat{B}_i}{\sum}\mu_m(t)m=0$ implique la nullit\'e des $\mu_m(t)$.

\noindent Il nous reste \`a montrer que tout $m\in \hat{A}$ est dans $F=(\underset{m\in \hat{B}_i}{\bigoplus}\ZZ[t,t^{-1}]E_i(m))\oplus (\underset{m\in \hat{A}-\hat{B}_i}{\bigoplus} \ZZ[t,t^{-1}]m)$. C'est clair si $m\in \hat{A}-\hat{B}_i$. Dans le cas $m\in \hat{B}_i$, montrons le par r\'ecurrence sur le $i$-poids de $m$. Si $m$ est de $i$-poids $0$, tous les $u_{i,a}(m)$ sont nuls et $m=E_i(m)$. Dans le cas g\'en\'eral, on a:
$$E_i(m)=m+f(m)=m+\underset{m'\in \hat{A}}{\sum}\lambda_{m'}(t)m'$$
avec $\lambda_{m'}(t)$ qui peut \^etre non nul seulement si le $i$-poids de $m'$ est strictement inf\'erieur \`a celui de $m$. Alors:
$$m=E_i(m)-\underset{m'\in \hat{A}-\hat{B}_i}{\sum}\lambda_{m'}(t)m'-\underset{m'\in \hat{B}_i}{\sum}\lambda_{m'}(t)m'$$
avec $E_i(m)\in F$, $\underset{m'\in \hat{A}-\hat{B}_i}{\sum}\lambda_{m'}(t)m'\in F$ et par hypoth\`ese de r\'ecurrence $\underset{m'\in \hat{B}_i}{\sum}\lambda_{m'}(t)m'\in F$.
\qed

\section{Les op\'erateurs $\hat{S_{t,i}}^l$}

\subsection{D\'efinition}

On consid\`ere les $\hat{\Yim}_t$-modules libres suivants ($i\in I$):
$$\hat{\Yim}^l_{t,i}=\underset{a\in\CC^*}{\bigoplus}\hat{\Yim}_t.S_{i,a}$$ 
On a alors une application naturelle: 
$$\tilde{\Pi}_{t,i}^l:\hat{\Yim}_{t,i}^l\rightarrow \Yim_{i}^l$$
d\'eduite de $\tilde{\Pi}_t$:
$$\tilde{\Pi}_{t,i}(\underset{a\in\CC^*}{\sum}\lambda_a.S_{i,a})=\underset{a\in\CC^*}{\sum}\tilde{\Pi}_t(\lambda_a).S_{i,a}$$

\begin{defi} On note $\hat{S}_{t,i}^l$ l'application $\ZZ[t,t^{-1}]$-lin\'eaire $\hat{S}_{t,i}^l:\hat{\Yim}_t\rightarrow \hat{\Yim}_{t,i}^l$ qui prend sur un mon\^ome $m\in \hat{A}$ la valeur:
$$\hat{S}_{t,i}^l(m)=m(\underset{a\in\CC^*/u_{i,a}(m)\geq 0}{\sum}(1+...+t^{2(u_{i,a}(m)-1)})S_{i,a}-\underset{a\in\CC^*/u_{i,a}(m)<0}{\sum}(t^{-2}+...+t^{2u_{i,a}(m)})S_{i,a})$$
\end{defi}
 
\begin{lem}
Le diagramme (1) suivant est commutatif:
$$\begin{array}{rcccl}
\hat{\Yim}_t&\stackrel{\hat{S}_{t,i}^l}{\longrightarrow}&\hat{\Yim}_{t,i}^l\\
\tilde{\Pi}_t\downarrow &&\downarrow&\tilde{\Pi}_{t,i}^l\\
\Yim&\stackrel{S_i^l}{\longrightarrow}&\Yim_{1,i}^l\\
\end{array}$$ 
\end{lem}

\demo

 Toutes les applications sont $\ZZ$-lin\'eaires, il suffit donc de regarder un mon\^ome $m\in\hat{A}$ et $\lambda(t)\in\ZZ[t,t^{-1}]$:
$$\begin{array}{ll}&(\tilde{\Pi}_{t,i}^l\circ \hat{S}^l_{t,i})(\lambda(t) m)
\\=&\lambda(1)\tilde{\Pi}_{t,i}^l(m(\underset{a\in\CC^*/u_{i,a}\geq 0}{\sum}(1+...+t^{2(u_{i,a}-1)})S_{i,a}-\underset{a\in\CC^*/u_{i,a}<0}{\sum}(t^{-2}+...+t^{2u_{i,a}})S_{i,a})
\\=&\lambda(1)\tilde{\Pi}_t(m)(\underset{a\in\CC^*}{\sum}u_{i,a}S_{i,a})
\\=&S_i^l(\lambda(1)\tilde{\Pi}_t(m))
\\=&S_i^l(\tilde{\Pi}_t(\lambda(t) m))\end{array}$$\qed

\subsection{Interpr\'etation de $\hat{S}_{t,i}^l$ en terme de d\'erivation}

\subsubsection{Des lois de bimodule sur $\hat{\Yim}^l_{t,i}$}

\begin{lem} Il existe sur $\hat{\Yim}^l_{t,i}$ une unique structure de bimodule pour la multiplication usuelle de $\hat{\Yim}_t$ telle que la structure \`a gauche soit la structure naturelle, que pour tout $a\in\CC^*$ et tout mon\^ome $m\in\hat{A}$:
$$S_{i,a}.m=t^{2u_{i,a}(m)}m.S_{i,a}\text{ , }S_{i,a}.t=t.S_{i,a}$$
\end{lem}

\demo

L'unicit\'e est claire, car la compatibilit\'e entre les structures \`a gauche et \`a droite impose, pour $\lambda_a\in\hat{\Yim}_t$, $\mu_m\in\ZZ[t,t^{-1}]$:
$$(\underset{a\in\CC^*}{\sum}\lambda_a.S_{i,a}).\underset{m}{\sum}\mu_mm=\underset{a\in\CC^*}{\sum}\lambda_a.(S_{i,a}.\underset{m}{\sum}\mu_mm)=\underset{a\in\CC^*}{\sum}\lambda_a.\underset{m}{\sum}\mu_mt^{2u_{i,a}(m)}m.S_{i,a}$$

\noindent Pour montrer que la structure $\ZZ[t,t^{-1}]$-lin\'eaire de module \`a droite est bien d\'efinie, il suffit de v\'erifier que pour deux mon\^omes $m_1,m_2\in\hat{A}$ on a $S_{i,a}.(m_1.m_2)=(S_{i,a}.m_1).m_2$. Ceci d\'ecoule du fait que $u_{i,a}(m_1.m_2)=u_{i,a}(m_1)+u_{i,a}(m_2)$. La compatibilit\'e entre les deux structures de modules est alors imm\'ediate.
\qed

\noindent On peut g\'en\'eraliser ce qui pr\'ec\`ede au cas d'une multiplication tordue sur $\hat{\Yim}_t$. Pour tout bicaract\`ere $d:\hat{A}\times\hat{A}\rightarrow\ZZ$, c'est-\`a-dire v\'erifiant:
$$d(m_1.m_2,m_3)=d(m_1,m_3)+d(m_2,m_3)\text{ , }d(m_1,m_2.m_3)=d(m_1,m_2)+d(m_1,m_3)$$
on a une loi de composition interne $*_d$ associative et $\ZZ[t,t^{-1}]$-lin\'eaire sur $\hat{\Yim}_t$ en posant:
$$m_1*_dm_2=t^{2d(m_1,m_2)}m_1.m_2$$
Remarquons que pour obtenir une loi associative, il suffit de demander que $d$ v\'erifie sur des mon\^omes $m_1,m_2,m_3$, la propri\'et\'e de cocycle:
$$-d(m_2,m_3)+d(m_1m_2,m_3)-d(m_1,m_2m_3)+d(m_1,m_2)=0$$
ce qui est le cas pour les bicaract\`eres.

\noindent On peut munir naturellement $\hat{\Yim}^l_{t,i}$ d'une structure de $\hat{\Yim}_t$-module \`a gauche pour la multiplication $*_d$ en posant pour $U\in\hat{\Yim}_t$:
$$U*_d(\underset{a\in\CC^*}{\sum}\lambda_a.S_a)=\underset{a\in\CC^*}{\sum}(U*_d\lambda_a).S_{i,a}$$
et de mani\`ere compl\`etement analogue au cas $d=0$ on fait de $\hat{\Yim}^l_{t,i}$ un bimodule en posant:
$$S_{i,a}*_dm=t^{2u_{i,a}(m)}m*_dS_{i,a}\text{ , }S_{i,a}*_dt=t*_dS_{i,a}$$

\subsubsection{$\hat{S}_{t,i}^l$ comme d\'erivation pour $*_d$}

\begin{prop} Pour $d$ bicarat\`ere, l'application $\hat{S}_{t,i}^l$ est une d\'erivation par rapport \`a la multiplication $*_d$:
$$\forall U,V\in \hat{\Yim}_t,\hat{S}_t^l(U*_dV)=U*_d\hat{S}_t^l(V)+\hat{S}_t^l(U)*_dV$$
\end{prop}

\demo

 Pour v\'erifier la propri\'et\'e de d\'erivation, et il suffit de montrer que pour deux mon\^omes $m,m'\in\hat{A}$ on a $\hat{S}^l_{t,i}(m*_dm')=\hat{S}^l_{t,i}(m)*_dm'+m*_d\hat{S}^l_{t,i}(m')$. Calculons en effet: 
$$\begin{array}{ll}&\hat{S}^l_{t,i}(m)*m'+m*\hat{S}^l_{t,i}(m')
\\=&m*_dm'
\\&(\underset{a\in\CC^{*}/u_{i,a}\geq 0}{\sum}(1+...+t^{2(u_{i,a}-1)})t^{2u_{i,a'}}S_{i,a}-\underset{a\in\CC^*/u_{i,a}<0}{\sum}(t^{-2}+...+t^{2u_{i,a}})t^{2u_{i,a'}}S_{i,a}
\\&+\underset{b\in\CC^{*}/u_{i,b}'\geq 0}{\sum}(1+...+t^{2(u_{i,b}'-1)})S_{i,b}-\underset{b\in\CC^*/u_{i,b}'<0}{\sum}(t^{-2}+...+t^{2u_{i,b}'})S_{i,b})
\\=&m*_dm'(\underset{a\in\CC^*/u_{i,a}+u_{i,a}'\geq 0}{\sum}(1+...+t^{2(u_{i,a}+u_{i,a}'-1)})S_{i,a}
\\&-\underset{a\in\CC^*/u_{i,a}+u_{i,a}'<0}{\sum}(t^{-2}+...+t^{2(u_{i,a}+u_{i,a}')})S_{i,a})
\\=&\hat{S}^l_{t,i}(m*_dm')
\end{array}$$
\qed

\noindent On a ainsi une caract\'erisation de $\hat{S}_{t,i}^l$ comme l'unique d\'erivation pour la loi usuelle, $\ZZ[t,t^{-1}]$-lin\'eaire, prenant les valeurs sur les g\'en\'erateurs:
$$\hat{S}_{t,i}^l(V_{i,a})=-t^{-2}V_{i,a}(S_{i,aq_i^{-1}}+S_{i,aq_i})\text{ et }\hat{S}_{t,i}^l(W_{j,a})=\delta_{i,j}W_{i,a}.S_{i,a}$$
$$\hat{S}_{t,i}^l(V_{j,a})=\delta_{C_{i,j},-1}V_{j,a}S_{i,a}+\delta_{C_{i,j},-2}V_{j,a}(S_{i,aq}+S_{i,aq^{-1}})+\delta_{C_{i,j},-3}V_{j,a}(S_{i,aq^{-2}}+S_{i,a}+S_{i,aq^2})$$
pour $j\neq i$.
\section{Les $t$-op\'erateurs d'\'ecrantage $\hat{S}_{t,i}$}

\subsection{D\'efinition de $\hat{S}_{t,i}$}

\begin{defi} On consid\`ere le $\ZZ[t,t^{-1}]$-sous module de $\hat{\Yim}^l_{t,i}$: 
$$\hat{F}_{t,i}=\underset{a\in\CC^*,m\in \hat{A}}{\sum}\ZZ[t,t^{-1}]m(V_{i,aq_i}t^{2u_{i,aq_i^2}(m)}S_{i,aq_i^2}-t^2S_{i,a})$$
Le module quotient obtenu est not\'e:
$$\hat{\Yim}_{t,i}=\hat{\Yim}_{t,i}^l/\hat{F}_{t,i}$$
L'application obtenue \`a partir de $\hat{S}^l_{t,i}$ par composition avec la projection $\hat{p}_{t,i}$ de $\hat{\Yim}_{t,i}^l$ sur $\hat{\Yim}_{t,i}$ est not\'ee $\hat{S}_{t,i}$.
\end{defi}

\noindent Nous allons montrer, en particulier dans le th\'eor\`eme \ref{justi}, que ces op\'erateurs peuvent \^etre consid\'er\'es comme des $t$-analogues des op\'erateurs d'\'ecrantage.

\begin{lem}
\noindent L'application $\tilde{\Pi}_{t,i}^l$ donne naturellement une application $\tilde{\Pi}_{t,i}$ rendant le diagramme (2) suivant commutatif:
$$\begin{array}{rcccl}
\hat{\Yim}_{t,i}^l&\longrightarrow&\hat{\Yim}_{1,i,t}\\
\tilde{\Pi}_{t,i}^l\downarrow &&\downarrow&\tilde{\Pi}_{t,i}\\
\Yim_{i}^l&\longrightarrow&\Yim_{i}\\
\end{array}$$
\end{lem}

\demo

Il suffit de v\'erifier que l'application $\ZZ$-lin\'eaire $\tilde{\Pi}_{t,i}^l$ passe au quotient. Pour 
$$x=\underset{a\in\CC^*,m\in \hat{A}}{\sum}\lambda_{m,a}m(V_{i,aq_i}t^{2u_{i,aq_i^2}(m)}S_{i,aq_i^2}-t^2S_{i,a})\in\hat{F}_{t,i}$$
on a:
$$\tilde{\Pi}_{t,i}(x)=\underset{a\in\CC^*}{\sum}\lambda_{m,a}(1)\Pi_t(m)(A_{i,aq_i}^{-1}.S_{i,aq_i^2}-S_{i,a})\in F_i$$
\qed

\begin{prop} On a le diagramme commutatif suivant:
$$\begin{array}{rcccl}
\hat{\Yim}_t&\stackrel{\hat{S}_{t,i}}{\longrightarrow}&\hat{\Yim}_{t,i}\\
\tilde{\Pi}_t\downarrow &&\downarrow&\tilde{\Pi}_{t,i}\\
\Yim&\stackrel{S_i}{\longrightarrow}&\Yim_{i}\\
\end{array}$$ 
\end{prop}

\demo 

La commutativit\'e du diagramme provient de la commutativit\'e des diagrammes (1) et (2). 
\qed

\noindent Soit $(\hat{\Yim}_t)_i=\ZZ[V_{i,a},W_{i,a}]_{a\in\CC^*}\subset\hat{\Yim}_t$. On d\'efinit alors $\pi_i:\hat{\Yim}\rightarrow (\hat{\Yim}_t)_i$ comme l'unique morphisme d'anneaux $\ZZ[t,t^{-1}]$-lin\'eaire tel que:
$$\pi_i(W_{i,a})=W_{i,a}\text{ , }\pi_i(V_{i,a})=V_{i,a}$$
$$\pi_i(W_{j,a})=1\text{ si $j\neq i$}$$
$$\pi_i(V_{j,a})=1\text{ si $C_{i,j}=0$}$$
$$\pi_i(V_{j,a})=W_{i,a}\text{ si $C_{i,j}=-1$}$$
$$\pi_i(V_{j,a})=W_{i,aq}W_{i,aq^{-1}}\text{ si $C_{i,j}=-2$}$$
$$\pi_i(V_{j,a})=W_{i,aq^2}W_{i,a}W_{i,aq^{-2}}\text{ si $C_{i,j}=-3$}$$
Pour un mon\^ome $m\in\hat{A}$, on a alors $u_{i,a}(m)=u_{i,a}(\pi_i(m))$ pour $a\in\CC^*$.

\begin{prop}\label{pasgene} Le noyau de l'application $\hat{S}_{t,i}$ contient $\hat{\mathfrak{K}}_{t,i}$.\end{prop}

\demo

Soit $m\in\hat{A}$ un mon\^ome $i$-dominant. Dans $S_{t,i}^l(E_i(m))$, on peut factoriser tous les termes par $m$, et notons $\frac{\hat{S}_{t,i}^l(E_i(m))}{m}\in \hat{\Yim}_{t,i}^l$ la quantit\'e obtenue. Elle ne d\'epend que des $u_{i,a}(m)$ ($a\in\CC^*$), et donc:
$$\frac{\hat{S}_{t,i}^l(E_i(m))}{m}=\frac{\hat{S}_{t,i}^l(E_i(\pi_i(m)))}{\pi_i(m)}$$
Remarquons de plus que les $u_{i,a}$ \'etant conserv\'es, on a: 
$$\hat{S}_{t,i}^l(E_i(m))\in \hat{F}_{t,i}\Leftrightarrow\hat{S}_{t,i}^l(E_i(\pi_i(m)))\in\hat{F}_{t,i}$$ 
En cons\'equence il nous suffit de montrer $\hat{S}_{t,i}(E_i(m))=0$ pour $m\in \hat{B}_i\cap (\hat{\Yim}_t)_i$. Mais alors tout se passe comme si on travaillait avec $\Glie=\U_{q_i}(sl_2)$. On est ainsi ramen\'e au cas $ADE$ qui sera \'etablie plus bas, ind\'ependamment de ce qui pr\'ec\`e, dans proposition \ref{passeu}.\qed

\subsection{Interpr\'etation de $\hat{S}_{t,i}$ comme d\'erivation dans le cas $ADE$}

Dans cette sous-partie on se restreint au cas o\`u $\Glie$ est de type $ADE$. On a alors tous les $q_i=q$ et la matrice de Cartan est sym\'etrique.

\subsubsection{Rappels \cite{Nab} et compl\'ements sur la loi $*$ de Nakajima}

On pose pour deux mon\^omes $m_1,m_2\in\hat{A}$:
$$\begin{array}{ll}d_N(m_1,m_2)&=\underset{i\in I,a\in\CC^*}{\sum}
(v_{i,aq}(m_1)u_{i,a}(m_2)+w_{i,aq}(m_1)v_{i,a}(m_2))
\\&=\underset{i\in I,a\in\CC^*}{\sum}
(u_{i,a}(m_1)v_{i,aq^{-1}}(m_2)+v_{i,a}(m_1)w_{i,aq^{-1}}(m_2))\end{array}$$
\noindent Ce bicaract\`ere, introduit par Nakajima dans \cite{Naa} et \cite{Nab}, permet comme pr\'ec\'edemment de d\'efinir une nouvelle multiplication sur 
$\hat{\Yim}_t$ en posant pour $m_1,m_2$ deux mon\^omes:
$$m_1*_{d_N}m_2=t^{2d_N(m_1,m_2)}m_1.m_2$$
avec $.$ la multiplication usuelle. On notera dans la suite simplement $d$ et $*$. Cette nouvelle multiplication $*$ n'est pas commutative.

\noindent Notons que $\hat{\mathfrak{K}}$ est une partie de $\hat{\Yim}_t$ stable pour la multiplication $*$ (\cite{Nab}).

\begin{lem}\label{utile} Soit $m\in \hat{A}$ un mon\^ome, $i\in I$ et $a\in\CC^*$. On a alors: 
$$V_{i,aq}*m=t^{2(u_{i,a}(m)-u_{i,aq^2}(m))}m*V_{i,aq}=t^{2u_{i,a}(m)}V_{i,aq}.m$$
\end{lem}

\noindent C'est une cons\'equence imm\'ediate de 
$$d(V_{i,aq},m)=u_{i,a}(m)\text{ et }d(m,V_{i,aq})=u_{i,aq^2}(m)$$

\begin{lem}\label{nakam} Soit $(m_1,...,m_p)$ des mon\^omes tels qu'il existe un $a\in\CC^*$ v\'erifiant pour tout $r$, $u_{i,a}(m_r)=1$ et $u_{i,b}(m_r)=0$ pour $b\neq a$. Alors il existe $\alpha\in\ZZ$ tel que:
$$(m_1*(1+V_{i,aq}))*(m_2*(1+V_{i,aq}))*...*(m_p*(1+V_{i,aq}))=t^{\alpha}m_1...m_p\underset{r=0..p}{\sum}t^{r(p-r)}\begin{bmatrix}p\\r\end{bmatrix}_tV_{i,aq}^r$$
\end{lem}

\demo

On proc\`ede par r\'ecurrence sur $p$ en s'appuyant sur le lemme \ref{utile}. Pour $p=1$, on a $m_1*V_{i,aq}=m_1V_{i,aq}$ et on a le r\'esultat avec $\alpha=1$. Ensuite dans le cas g\'en\'eral:
$$\begin{array}{ll}&(m_1*(1+V_{i,aq}))*(m_2*(1+V_{i,aq}))*...*(m_{p+1}*(1+V_{i,aq}))
\\&=t^{\alpha}(m_1+m_1V_{i,aq})*(m_2...m_{p+1}\underset{r=0..p}{\sum}t^{r(p-r)}\begin{bmatrix}p\\r\end{bmatrix}_tV_{i,aq}^r)
\\&=t^{\alpha+2d(m_1,m_2...m_{p+1})}m_1m_2...m_{p+1}\underset{r=0..p}{\sum}t^{r(p-r)}\begin{bmatrix}p\\r\end{bmatrix}_tV_{i,aq}^r
\\&+t^{\alpha+2d(m_1,m_2...m_{p+1})}m_1m_2...m_{p+1}\underset{r=0..p}{\sum}t^{r(p-r)}\begin{bmatrix}p\\r\end{bmatrix}_tt^{2p-2r}V_{i,aq}^{r+1}
\\&=t^{\alpha+2d(m_1,m_2...m_{p+1})}m_1m_2...m_{p+1}\underset{r=0..p+1}{\sum}(t^{r(p-r)}\begin{bmatrix}p\\r\end{bmatrix}_t+t^{(r-1)(p-r+1)}\begin{bmatrix}p\\r-1\end{bmatrix}_tt^{2p-2r+2})V_{i,aq}^r
\end{array}$$
Et on conclut en remarquant:
$$t^{r(p-r)}\begin{bmatrix}p\\r\end{bmatrix}_t+t^{(r-1)(p-r+1)}\begin{bmatrix}p\\r-1\end{bmatrix}_tt^{2p-2r+2}=t^{r(p+1-r)}\begin{bmatrix}p+1\\r\end{bmatrix}_t$$
\qed

On a en particulier le r\'esultat:

\begin{lem}\label{nak} Soit $m$ un mon\^ome tel qu'il existe un $a\in\CC^*$ v\'erifiant $u_{i,a}(m)=1$ et $u_{i,b}(m)=0$ pour $b\neq a$. Alors pour $l\geq 0$:
$$[m(1+V_{i,aq})]^{*l}=m^l\underset{r=0..l}{\sum}t^{r(l-r)}\begin{bmatrix}l\\r\end{bmatrix}_tV_{aq}^r$$
\end{lem}

On peut exprimer les $E_i(m)$ en utilisant la loi $*$:

\begin{prop}\label{spee} On fixe un $i\in I$. Soit $m\in\hat{B}_i$ un mon\^ome $i$-dominant. Pour $a\in\CC^*$, on consid\`ere la suite $(Z_{i,a})=(Z_{i,a,l})_{1\leq l\leq z_{i,a}}$ form\'ee de 
$$z_{i,a}=w_{i,a}+\underset{j/C_{j,i}=-1}{\sum}v_{j,a}=u_{i,a}+v_{i,aq}+v_{i,aq^{-1}}$$ 
termes o\`u $W_{i,a}$ appara\^it $w_{i,a}$ fois et pour $j$ tel que $C_{j,i}=-1$, $V_{j,a}$ appara\^it $v_{j,a}$ fois:  
$$\{W_{i,a},...,W_{i,a},V_{j_1,a},...,V_{j_1,a},V_{j_2,a},...,V_{j_m,a}\}=\{Z_{i,a,1},...,Z_{i,a,z_{i,a}}\}$$
Alors il existe un unique $\beta\in\ZZ$ tel que:
$$t^{\beta}E_i(m)=(\underset{j\neq i,a\in\CC^*}{\prod^*}W_{j,a})
*(\underset{j/C_{j,i}=0,a\in(\CC^*/q^{2\ZZ}),r\in\ZZ}{\overset{\rightarrow}{\prod^*}}V_{j,aq^{2r}})
*(\underset{a\in(\CC^*/q^{2\ZZ})}{\prod^*}\underset{r\in\ZZ}{\overset{\rightarrow}{\prod^*}}m_{i,a,r})$$
avec:
$$\begin{array}{ll}m_{i,a,r}
&=(Z_{i,a,1}*(1+V_{i,aq}))*...*(Z_{i,a,u_{i,a}}*(1+V_{i,aq}))
\\&*(Z_{i,a,u_{i,a}+1}*V_{i,aq}*Z_{i,aq^2,u_{i,aq^2}+v_{i,aq^3}+1})*...*(Z_{i,a,u_{i,a}+v_{i,aq}}*V_{i,aq}*Z_{i,aq^2,u_{i,aq^2}+v_{i,aq^3}+v_{i,aq}})\end{array}$$
\end{prop}

\demo

On commence par expliciter $E_i(m)$:
$$E_i(m)=m\underset{a\in\CC^*}{\prod}(\underset{r_a=0..u_{i,a}(m)}{\sum}t^{r_a(u_{i,a}(m)-r_a)}\begin{bmatrix}u_{i,a}(m)\\r_a\end{bmatrix}_tV_{i,aq}^{r_a})$$
Si on ne tient pas compte des $t$, l'expression annonc\'ee est correcte puisqu'on a le bon nombre $v_{i,aq}$ de $V_{i,aq}$ et tous les $Z_{i,a,l}$ pour $l=1...u_{i,a}+v_{i,aq}+v_{i,aq^{-1}}$. Le seul probl\`eme est l'inhomog\`en\'e\"it\'e de $E_i(m)$ du fait des puissances de $V_{i,aq}$.

\noindent Les seuls facteurs de $m$ qui contribuent aux $u_{i,b}(m)$ sont les $W_{i,a},V_{i,a}$ et les $V_{j,a}$ avec $C_{j,i}=-1$. Mais ce sont exactement les facteurs qui posent probl\`eme avec $V_{i,aq}$ d'apr\`es le lemme \ref{utile}. On en d\'eduit une premi\`ere expression:
$$t^{\alpha}E_i(m)=(\underset{j\neq i,a\in\CC^*}{\prod^*}W_{j,a})
*(\underset{j/C_{j,i}=0,a\in(\CC^*/q^{2\ZZ}),r\in\ZZ}{\overset{\rightarrow}{\prod^*}}V_{j,aq^{2r}})
*(m'\underset{r_a=0..u_{i,a}(m)}{\sum}t^{r_a(u_{i,a}(m)-r_a)}\begin{bmatrix}u_{i,a}(m)\\r_a\end{bmatrix}_tV_{i,aq}^{r_a})$$
avec 
$$\begin{array}{ll}m'&=(\underset{a\in\CC^*}{\prod}W_{i,a}^{w_{i,a}}V_{i,a}^{v_{i,a}})(\underset{a\in\CC^*,j/C_{i,j}=-1}{\prod}V_{j,a}^{v_{j,a}})
\\&=\underset{a\in \CC^*}{\prod}(\underset{l=1..u_{i,a}}{\prod}Z_{i,a,l})(\underset{r=1...v_{i,aq}}{\prod}Z_{i,a,u_{i,a}+r}V_{i,aq}Z_{i,aq^2,u_{i,aq^2}+v_{i,aq^3}+r})\end{array}$$
Il nous suffit donc de montrer qu'il existe $\gamma\in\ZZ$ tel que:
$$E_i(m')=m'\underset{r_a=0..u_{i,a}(m)}{\sum}t^{r_a(u_{i,a}(m)-r_a)}\begin{bmatrix}u_{i,a}(m)\\r_a\end{bmatrix}_tV_{i,aq}^{r_a}=t^{-\gamma}(\underset{a\in\CC^*/q^{2\ZZ}}{\prod^*}\underset{r\in\ZZ}{\overset{\rightarrow}{\prod^*}}m_{i,a,r})$$
Or d'apr\`es le lemme \ref{nakam}, le facteur $Z_{i,a,1}...Z_{i,a,u_{i,a}(m)}\underset{r_a=0..u_{i,a}(m)}{\sum}t^{r_a(u_{i,a}(m)-r_a)}\begin{bmatrix}u_{i,a}(m)\\r_a\end{bmatrix}_tV_{i,aq}^{r_a}$ est \'egal \`a une puissance de $t$ pr\`es \`a $(Z_{i,a,1}*(1+V_{i,aq}))*...*(Z_{i,a,u_{i,a}}*(1+V_{i,aq}))$. Il ne reste plus qu'\`a v\'erifier que les facteurs restant $\underset{a\in \CC^*}{\prod}(\underset{r=1...v_{i,aq}}{\prod}Z_{i,a,u_{i,a}+r}V_{i,aq}Z_{i,aq^2,u_{i,aq^2}+v_{i,aq^3}+r})$ ne posent pas de probl\`eme vis \`a vis de l'inhomog\`en\'e\"it\'e en puissances de $V_{i,aq}$, mais c'est le cas car pour tout $r\in\ZZ$: 
$$u_{i,aq^{2r}}(Z_{i,a,u_{i,a}+l}*V_{i,aq}*Z_{i,aq^2,u_{i,aq^2}+v_{i,aq^3}+l})=0$$
\qed

\subsubsection{Une structure de bimodule sur $\hat{\Yim}_{t,i}$ pour la loi $*$}

\begin{lem} Le sous $\ZZ[t,t^{-1}]$ module $\hat{F}_{t,i}$ de $\hat{\Yim}^l_{t,i}$ est en fait un sous-module \`a gauche pour la loi $*$:
$$\hat{F}_{t,i}=\underset{a\in\CC^*}{\sum}\hat{\Yim}_t*(V_{i,aq}.S_{i,aq^2}-t^2S_{i,a})$$
et m\^eme un sous-bimodule $\hat{\Yim}_t*\hat{F}_{t,i}=\hat{F}_{t,i}*\hat{\Yim}=\hat{F}_{t,i}$.
\end{lem}

\demo

La premi\`ere propri\'et\'e d\'ecoule directement du lemme \ref{utile} qui donne pour $m\in\hat{A}$:
$$m*(V_{i,aq}.S_{i,aq^2}-t^2S_{i,a})=t^{2u_{i,aq_i^2}(m)}mV_{i,aq_i}S_{i,aq_i^2}-t^2mS_{i,a}$$
Pour la propri\'et\'e de sous-bimodule, soit $m\in \hat{A}$ un mon\^ome. En utilisant le lemme \ref{utile}, on a pour $\lambda_a\in \hat{\Yim}_t$:
$$\begin{array}{ll}\lambda_a*(V_{i,aq}.S_{i,aq^2}-t^2S_{i,a})*m
&=\lambda_a*(t^{2u_{i,aq^2}(m)}V_{i,aq}*m.S_{i,aq^2}-t^{2+2u_{i,a}(m)}m.S_{i,a})
\\&=t^{2u_{i,a}(m)}\lambda_a*m*(V_{i,aq}.S_{i,aq^2}-t^2S_{i,a})\in \hat{F}_{t,i}\end{array}$$
\qed

\noindent On peut ainsi munir naturellement $\hat{\Yim}_{t,i}$ d'une structure de bimodule.

\subsubsection{$\hat{S}_{t,i}$ est une d\'erivation}

\noindent Le r\'esultat suivant, qui justifie entre autre les constructions pr\'ec\'edentes, permet en particulier d'obtenir la proposition \ref{pasgene}:

\begin{prop}\label{passeu} L'application $\hat{S}_{t,i}$ est une d\'erivation pour le produit $*$ et son noyau contient $\hat{\mathfrak{K}}_{t,i}$.\end{prop}

\demo

La propri\'et\'e de d\'erivation est conserv\'ee: en effet pour $U,V\in\hat{\Yim}_t$ on a:
$$\hat{S}_{t,i}(U*V)
=\hat{p}_{t,i}(U*\hat{S}_{t,i}^l(V))+\hat{p}_{t,i}(\hat{S}_{t,i}(U)*V)
=U*\hat{p}_{t,i}(\hat{S}_{t,i}^l(V))+\hat{p}_{t,i}(\hat{S}_{t,i}(U))*V
=U*\hat{S}_{t,i}(V)+\hat{S}_{t,i}(U)*V$$
Pour montrer que $\hat{\mathfrak{K}}_{t,i}\subset \text{Ker}(\hat{S}_{t,i})$, consid\'erons un mon\^ome dominant $m$, et d\'ecomposons en utilisant la proposition \ref{spee} sous la forme d'un produit pour $*$. En utilisant la propri\'et\'e de d\'erivation de $\hat{S}_{t,i}$, il nous suffit d'obtenir que chacun des termes est annul\'e. Or pour $a\in\CC^*$:
$$\hat{S}_{t,i}(Z_{i,a,k}*(1+V_{i,aq}))=Z_{i,a,k}.S_{i,a}-t^{-2}Z_{i,a,k}*V_{i,aq}.S_{i,aq^2}=Z_{i,a,k}*(S_{i,a}-t^{-2}V_{i,aq}S_{i,aq^2})=0$$
$$\hat{S}_{t,i}^l(Z_{i,a,k}*V_{aq}*Z_{i,aq^2,k'})=0$$
car pour tout $b\in\CC^*$, $u_{i,b}(Z_{i,a,k}*V_{aq}*Z_{i,aq^2,k'})=0$\qed

\subsection{Interpr\'etation de $\hat{S}_{t,i}$ dans le cas g\'en\'eral}

Dans le cas g\'en\'eral, on ne dispose pas de bicarat\`ere v\'erifiant les deux relations fondamentales du cas $ADE$ pour tout $i\in I$
$$d(V_{i,aq_i},m)=u_{i,a}(m)\text{ et }d(m,V_{i,aq_i})=u_{i,aq_i^2}(m)$$
Par exemple, pour $\Glie$ de type $B_2$, on a $C=\begin{pmatrix} 2 & -1 \\ -2 & 2 \end{pmatrix}$ (la matrice de Cartan dans \cite{Fre} est la transpos\'ee de celle de \cite{bou}), $q_1=q^2$, $q_2=q$ et:
$$0=u_{1,aq^{-2}}(V_{2,a})\neq u_{2,aq}(V_{1,a})=1$$
On ne peut pas traiter tous les op\'erateurs simultan\'ement, mais on peut cependant les interpr\'eter individuellement en posant pour chaque $i\in I$:
$$d_i(m_1,m_2)=\underset{a\in\CC^*}{\sum}(v_{i,aq_i}(m_1)u_{i,a}(m_2)+w_{i,aq_i}(m_1)v_{i,a}(m_2))+\underset{a\in\CC^*}{\sum}(u_{i,aq_i}-w_{i,aq_i}+v_{i,a}+v_{i,aq_i^2})(m_1)v_{i,a}(m_2)$$
$$=\underset{a\in\CC^*}{\sum}(u_{i,aq_i}(m_1)v_{i,a}(m_2)+v_{i,aq_i}(m_1)w_{i,a}(m_2))+\underset{a\in\CC^*}{\sum}v_{i,aq_i}(m_1)(u_{i,a}-w_{i,a}+v_{i,aq_i^{-1}}+v_{i,aq_i})(m_2)$$
Il d\'ecoule alors de la d\'efinition:
\begin{lem} Pour tout $m\in \hat{A}$:
$$d_i(V_{i,aq_i},m)=u_{i,a}(m)\text{ et }d_i(m,V_{i,aq_i})=u_{i,aq_i^2}(m)$$
\end{lem}
\noindent On note $*_i$ la loi sur $\hat{\Yim}_t$ associ\'ee au bicarat\`ere $d_i$. On montre alors de la m\^eme mani\`ere que dans le cas $ADE$:
\begin{prop} Le sous $\ZZ[t,t^{-1}]$ module $\hat{F}_{t,i}$ de $\hat{\Yim}^l_{t,i}$ est en fait un sous-module \`a gauche pour la loi $*_i$:
$$\hat{F}_{t,i}=\underset{a\in\CC^*}{\sum}\hat{\Yim}_t*_i(V_{i,aq}.S_{i,aq^2}-t^2S_{i,a})$$
et m\^eme un sous-bimodule $\hat{\Yim}_t*_i\hat{F}_{t,i}=\hat{F}_{t,i}*_i\hat{\Yim}=\hat{F}_{t,i}$.

\noindent L'application $\hat{S}_{t,i}$ est une d\'erivation pour le produit $*_i$.
\end{prop}

\subsection{D\'emonstration du th\'eor\`eme \ref{justi}}

On retourne au cas g\'en\'eral pour $\Glie$ simple quelconque.

\begin{thm}\label{justi} On a $\hat{\mathfrak{K}}_{t,i}= \text{Ker}(\hat{S}_{t,i})$.\end{thm}

\demo

La premi\`ere inclusion $\hat{\mathfrak{K}}_{t,i}\subset \text{Ker}(\hat{S}_{t,i})$ est d\'ej\`a connue dans la proposition \ref{pasgene}.

\noindent Supposons par l'absurde qu'on n'ait pas \'egalit\'e. Alors on consid\`ere un $x\in \text{Ker}(\hat{S}_{t,i})-\hat{\mathfrak{K}}_{t,i}$ qu'on d\'ecompose en utilisant le lemme \ref{sdirecte} sur $\hat{\Yim}_t=\hat{\mathfrak{K}}_{t,i}\oplus \underset{m\in \hat{A}-\hat{B}_i}{\bigoplus} \ZZ[t,t^{-1}]m$ sous la forme $x=v+u$ avec $u\neq 0$. On note l'\'ecriture de $u$:
$$u=\underset{m\in M}{\sum}\lambda_mm$$
avec $M\subset \hat{A}-\hat{B}_i$, $\lambda_m\in\ZZ[t,t^{-1}]$ et $\lambda_m\neq 0$ pour $m\in M$.
Alors $x,v\in \text{Ker}(\hat{S}_{t,i})$, donc $u$ est un \'el\'ement non nul de $\underset{m\in \hat{A}-\hat{B}_i}{\bigoplus} \ZZ[t,t^{-1}]m\cap \text{Ker}(\hat{S}_{t,i})$. Pour $m\in \hat{A}-\hat{B}_i$, notons $N_m$ le nombre de classe $R\in \CC^*/q^{2\ZZ}$ tel qu'il existe $a\in R$ v\'erifiant $u_{i,a}(m)<0$. Tous les mon\^omes $m$ de $\hat{A}-\hat{B}_i$ v\'erifient $N_m\geq 1$. Soit $m_0\in M$ avec $N_{m_0}$ minimal parmi les $N_{m}$ pour $m\in M$. Soit alors $a\in\CC^*$ tel que $u_{i,a}(m_0)< 0$ et pour $r<0$, $u_{i,aq^{2r}}(m_0)\geq 0$. Lorsqu'on calcule
$$\hat{S}_{t,i}(u)=0=\underset{m\in M}{\sum}\lambda_mm(\underset{b\in\CC^*/u_{i,b}(m)\geq 0}{\sum}(1+t^2+...+t^{2(u_{i,b}(m)-1)})S_{i,b}-\underset{b\in\CC^*/u_{i,b}(m)<0}{\sum}(t^{-2}+...+t^{2u_{i,b}(m)})S_{i,b})$$
on voit appa\^itre le terme $-\lambda_{m_0}m_0(t^{-2}+...+t^{2u_{i,a}(m)})S_{i,a}$. Ce terme doit \^etre annul\'e par projection sur $\hat{\Yim}_{t,i}$. Les termes qui vont l'annuler peuvent provenir soit d'un $S_{i,aq_i^{2r}}$ avec $r<0$, soit d'un $S_{i,aq_i^{2r}}$ avec $r>0$. Dans le premier cas on a un mon\^ome $m_1\in M$ tel que $m_1V_{i,aq_i^{-1}}V_{i,aq_i^{-3}}...V_{i,aq_i^{2r+1}}=m_0$, dans le deuxi\`eme on a un mom\^ome $m_1=m_0V_{i,aq_i}...V_{i,aq_i^{2r-1}}\in M$. On peut ainsi d\'efinir une suite de mon\^omes $m_p$ tant que $u_{i,a}(m_p)<0$. Les termes de la suite sont distincts deux \`a deux, car \`a chaque op\'eration soit on ajoute des $V_{i,aq_i^{2r+1}}$ avec $r\geq 0$, soit on enl\`eve des $V_{i,aq_i^{2r+1}}$ avec $r<0$. Notons aussi qu'\`a chaque op\'eration on ne diminue pas les $u_{i,aq_i^{2r}}$ avec $r<0$, et on n'augmente pas $N$. Comme $M$ est fini, la suite se termine sur un $m_P\in M$ qui v\'erifie $u_{i,aq_i^{2r}}(m_P)\geq 0$ pour $r\leq 0$, et $N_{m_P}=N_{m_0}$. En notant $m^0=m_0$ et $m^1=m_P$, ce nouveau proc\'ed\'e donne une suite $m^j$ telle que $\text{min}\{r/u_{i,aq^{2r}}<0\}$ est strictement croissante. Par finitude de $M$, la suite se termine sur un $m^{P'}\in M$ tel que $u_{i,aq^{2r}}\geq 0$ pour tout $r\in \ZZ$ et les autres classes de $\CC^*/q_i^{2\ZZ}$ n'ont pas \'et\'e modifi\'ees. Donc $N_{m_{P'}}<N_{m_0}$, contradiction.\qed

\section{Op\'erateurs d'\'ecrantage pour l'anneau $\Yim_t$}

\subsection{Rappels et compl\'ements}

\subsubsection{L'anneau $\Yim_t$}

En suivant Nakajima \cite{Naa}, \cite{Nab} on consid\`ere l'anneau ``interm\'ediaire'' entre $\hat{\Yim}_t$ et $\Yim$:
$$\Yim_t=\ZZ[t,t^{-1},Y_{i,a},Y_{i,a}^{-1}]_{i\in I,a\in\CC^*}$$
On a un morphisme d'anneaux canonique:
$$\Pi_t:\Yim_t\rightarrow\Yim$$
$$Y_{i,a}^{\pm}\mapsto Y_{i,a}^{\pm} \text{ et } t \mapsto 1$$
Pour passer de $\hat{\Yim}_t$ \`a $\Yim$, on peut consid\'erer pour tout bicaract\`ere $d$ l'application $\hat{\Pi}_d:\hat{\Yim}_t\rightarrow \Yim_t$ qui est $\ZZ[t,t^{-1}]$-lin\'eaire, et qui v\'erifie:
$$m=\underset{i\in I,a\in\CC^*}{\prod}V_{i,a}^{v_{i,a}(m)}W_{i,a}^{w_{i,a}(m)}\mapsto t^{-d(m,m)}\underset{i\in I,a\in\CC^*}{\prod}Y_{i,a}^{u_{i,a}(m)}\text{ et }t\mapsto 1$$
On a toujours $\tilde{\Pi}_t=\Pi_t\circ \hat{\Pi}_d$. 

\noindent Dans le cas du bicaract\`ere trivial $d=0$, on note $\hat{\Pi}_0=\hat{\Pi}_t$ et c'est alors un morphisme d'anneaux. Dans le cas $ADE$, on peut prendre $d_N$ et on retrouve l'application $\hat{\Pi}=\hat{\Pi}_{d_N}$ de \cite{Nab}.

\begin{lem}\label{ordre} Un produit $p=\underset{i\in I,a\in\CC^*}{\prod}A_{i,a}^{v_{i,a}}\in\Yim$ avec les $v_{i,a}\in\ZZ$ est \'egal \`a $1$ si et seulement si tous les $v_{i,a}$ sont nuls.

\noindent En cons\'equence on d\'efinit une relation d'ordre partiel sur l'ensemble $A$ des mon\^omes de $\Yim$ en posant
:
$$m\leq m'\Leftrightarrow m'/m \text{ est un mon\^ome en }A_{i,a}^{-1}$$
\end{lem}

\demo

Supposons par l'absurde qu'un tel produit $p$ peut \^etre \'egal \`a $1$ avec des $v_{i,a}\neq 0$. Consid\'erons alors un $a$ tel qu'il existe un $i\in I$ avec $v_{i,a}\neq 0$ mais pour $m\in \ZZ$ stritement positif, pour $j\in I$, $v_{j, aq^m}=0$. Parmi ces $i$, on en choisit un tel que la longueur de la racine associ\'ee soit maximale. Dans $p$, le facteur $Y_{i,aq_i}^{-v_{i,a}}$ doit se simplifier avec un autre facteur. Cependant par d\'efinition de $a$ il ne peut pas venir de $A_{i,aq_i^2}^{v_{i,aq_i^2}}$. Il reste donc les possibilit\'es suivantes:

il provient d'un $A_{j,aq_i}^{v_{j,aq_i}}$ avec $C_{i,j}=-1$, $j\neq i$. Alors $v_{j,aq_i}\neq 0$, contradiction.

il provient d'un $A_{j,aq_iq^{-1}}^{v_{j,aq_iq^{-1}}}$ avec $C_{i,j}=-2$, $j\neq i$, ce qui impose $v_{j,aq_iq^{-1}}\neq 0$. Comme $C_{i,j}=-2$, les racines associ\'ees \`a $i$ et $j$ ne sont pas de m\^eme longueur, et donc en utilisant l'hypoth\`ese sur $i$, on a $r_i>r_j\geq 1$. Alors $q_iq^{-1}=q\text{ ou }q^2$, donc $v_{j,aq}\neq 0$ ou $v_{j,aq^2}\neq 0$, ce qui n'est pas possible d'apr\`es le choix de $i$.

il provient d'un $A_{j,aq_iq^{-2}}^{v_{j,aq_iq^{-2}}}$ avec $C_{i,j}=-3$, $j\neq i$, ce qui impose $v_{j,aq_iq^{-2}}\neq 0$. On est dans le cas o\`u $\Glie$ est de type $G_2$. Les racines associ\'ees \`a $i$ et $j$ ne sont pas de m\^eme longueur, et $r_i/r_j=3\text{ ou }\frac{1}{3}$. Si $r_i=3$, on a $v_{j,aq}\neq 0$ ce qui est contraire au choix de $i$. Si $r_i=1$, on a $v_{j,aq^{-1}}\neq 0$ et $v_{j,aq^m}=0$ pour $m\geq 0$. Mais alors on ne peut pas annuler $Y_{j,aq^{-1}q_j}^{v_{j,aq^{-1}}}=Y_{j,aq^2}^{v_{j,aq^{-1}}}$.

\noindent Pour que $\leq$ soit bien une relation d'ordre, la propri\'et\'e la moins \'evidente est l'antisym\'etrie qui est assur\'ee par ce qui pr\'ec\`ede.\qed

\subsubsection{Quelques notations}

On note $A$ l'ensemble des mon\^omes de $\Yim$, $B_i\subset A$ l'ensemble des mon\^omes $i$-dominants de $\Yim$. Pour $\underset{a\in\CC^*,j}{\prod}Y_{j,a}^{u_{j,a}}=m\in B_i$, on pose:
$$E_{0,i}(m)=m\underset{a\in\CC^*}{\prod}\underset{r_a=0..u_{i,a}}{\sum}t^{r_a(u_{i,a}-r_a)}\begin{bmatrix}u_{i,a}\\r_a\end{bmatrix}_tA_{i,aq_i}^{-r_a}$$
Remarquer que si de plus $m\in B_i\cap \hat{\Pi}_t(\hat{B}_i)=B_i'$, on a, en posant $m=\hat{\Pi}_t(m')$, l'\'egalit\'e $E_{0,i}(m)=\hat{\Pi}_t(E_i(m'))$. 

\noindent On note $\mathfrak{K}_i$ le $\ZZ[t,t^{-1}]$-module engendr\'e par les $E_{0,i}(m)$ avec $m\in B_i$. On a $\hat{\Pi}_t(\hat{\mathfrak{K}}_i)\subset\mathfrak{K}_i$ mais on n'a pas \'egalit\'e dans le cas g\'en\'eral.

\noindent Pour $i\in I$, on obtient de la m\^eme mani\`ere que dans le lemme \ref{sdirecte} une d\'ecomposition en somme directe de $\ZZ[t,t^{-1}]$-modules:
\begin{lem}\label{sdirectesimp}
$$\Yim_t=\mathfrak{K}_{t,i}\oplus \underset{m\in A-B_i}{\bigoplus} \ZZ[t,t^{-1}]m=(\underset{m\in B_i}{\bigoplus}\ZZ[t,t^{-1}]E_{i,0}(m))\oplus (\underset{m\in A-B_i}{\bigoplus} \ZZ[t,t^{-1}]m)$$
\end{lem}

\begin{lem}\label{mathfra} On a l'\'egalit\'e:
$$\hat{\Pi}_t(\hat{\mathfrak{K}_t})=\underset{i\in I}{\bigcap}\mathfrak{K}_{t,i}$$
et on notera $\mathfrak{K}_t$ cette sous-partie de $\Yim_t$.
\end{lem}

\demo

On sait d\'ej\`a:
$$\hat{\Pi}_t(\hat{\mathfrak{K}_t})=\hat{\Pi}_t(\underset{i\in I}{\bigcap}\hat{\mathfrak{K}}_i)\subset \underset{i\in I}{\bigcap}\hat{\Pi}_t(\hat{\mathfrak{K}}_i)\subset \underset{i\in I}{\bigcap}\mathfrak{K}_i$$
Consid\'erons \`a pr\'esent $x\in \underset{i\in I}{\bigcap}\mathfrak{K}_i$. Soit $m\in A$ un mon\^ome maximal parmi ceux qui interviennent dans $x$ pour la relation d'ordre $\leq$ du lemme \ref{ordre}. Pour chaque $i\in I$, $m$ provient d'un certain $E_{0,i}(m')$ avec $m'$ $i$-dominant, ce qui impose $m\leq m'$. On a donc $m=m'$ et $m$ est $i$-dominant pour tout $i\in I$. Il est donc de la forme: 
$$m=\underset{i\in I,a\in \CC^*}{\prod}Y_{i,a}^{u_{i,a}(m)}=\hat{\Pi}_t(\underset{i\in I,a\in \CC^*}{\prod}W_{i,a}^{u_{i,a}(m)})\in\hat{\Pi}_t(\hat{\Yim}_t)$$
car les $u_{i,a}(m)\geq 0$. Si on suppose $m\neq 1$ (soit $x\not\in\ZZ[t,t^{-1}]$) et on consid\`ere $i_0\in I$ tel que $wt_{i_0}(m)\neq 0$, on a dans l'\'ecriture de $x$ dans $\mathfrak{K}_{i_0}$ le mon\^ome $m$ qui ne peut provenir que de 
$$E_{0,{i_0}}(m)=\hat{\Pi}_t(E_{i_0}(\underset{i\in I,a\in\CC^*}{\prod}W_{i,a}^{u_{i,a}(m)}))\in \hat{\Pi}_t(\hat{\mathfrak{K}_t})$$
On peut alors enlever de $x$ le terme $E_{0,{i_0}}(m)$ avec son coefficient de $\ZZ[t,t^{-1}]$, et on se ram\`ene \`a un \'el\'ement de $\underset{i\in I}{\bigcap}\mathfrak{K}_{t,i}$ faisant intervenir strictement moins de mon\^ome, ce qui permet de conclure par r\'ecurrence.
\qed

\noindent Pour d\'efinir les $E_{0,i}'(m)$ analogues des $E_{0,i}(m)$ relatifs \`a $\hat{\Pi}$, on consid\`ere pour $i\in I$ en suivant \cite{Naa}:
$$\phi_i:\mathfrak{K}_{t,i}\rightarrow \Yim_t$$
d\'efinie comme l'application $\ZZ[t,t^{-1}]$-lin\'eaire telle que pour $m\in B_i$:
$$E_{0,i}(m)=\sum\lambda_M(t)M\mapsto \sum\lambda_M(t)t^{-\alpha(m,M)}M$$
avec pour $M=m\underset{a\in\CC^*}{\prod}A_{i,a}^{-r_a}$ ($r_a\geq 0$) qui intervient effectivement dans $E_{0,i}(m)$:
$$\alpha(m,M)=\underset{a\in\CC^*}{\sum}r_a(u_{i,aq_i^{-1}}(m)+u_{i,aq_i}(m)-r_a-r_{aq_i^{-2}})$$
On pose alors pour $m\in B_i$:
$$E_{0,i}'(m)=\phi_i(E_{0,i}(m))$$
Cette d\'efinition est motiv\'ee par le lemme:
\begin{lem} Soit, dans le cas $ADE$, $m\in B_i'$. On a, si $m=\hat{\Pi}_t(m')$, l'\'egalit\'e: 
$$E_{0,i}'(m)=t^{-d(m',m')}\hat{\Pi}(E_i(m'))$$
\end{lem}

\demo

Il suffit de calculer en notant $E_i(m')=\sum\lambda_M(t)M$ avec $M=m'\underset{a\in\CC^*}{\prod}V_{i,a}^{r_{a,M}}$:
$$\hat{\Pi}(E_i(m'))
=\hat{\Pi}(\sum\lambda_{M}(t)M)
=\sum \lambda_{M}(t) t^{-d(M,M)}m\underset{a\in\CC^*}{\prod}A_{i,a}^{-r_{a,M}}$$
puis:
$$d(M,M)=d(m',m')+\underset{a\in\CC^*}{\sum}r_{a,M}(d(V_{i,a},m')+d(m',V_{i,a})+d(V_{i,a},M/m'))$$
$$=d(m',m')+\underset{a\in\CC^*}{\sum}r_{a,M}(u_{i,aq_i^{-1}}(m')+u_{i,aq_i}(m')-r_{a,M}-r_{aq_i^{-2},M})=d(m',m')+\alpha(m,\hat{\Pi}_t(M))$$
ce qui donne
$$\hat{\Pi}(E_i(m'))=t^{-d(m',m')}\phi_i(\hat{\Pi}_t(E_i(m')))=t^{-d(m',m')}\phi_i(E_{0,i}(m))$$
\qed

\noindent On note alors $\mathfrak{K}_i'$ le $\ZZ[t,t^{-1}]$-module engendr\'e par les $E_{0,i}(m)'$ avec $m\in B_i$. 

\noindent Pour $i\in I$, on obtient de la m\^eme mani\`ere que dans le lemme \ref{sdirecte} une d\'ecomposition en somme directe de $\ZZ[t,t^{-1}]$-modules:
$$\Yim_t=\mathfrak{K}_{t,i}'\oplus \underset{m\in A-B_i}{\bigoplus} \ZZ[t,t^{-1}]m=(\underset{m\in B_i}{\bigoplus}\ZZ[t,t^{-1}]E_{i,0}'(m))\oplus (\underset{m\in A-B_i}{\bigoplus} \ZZ[t,t^{-1}]m)$$

\noindent Dans le cas $ADE$, on a $\hat{\Pi}(\hat{\mathfrak{K}}_i)\subset\mathfrak{K}_i'$ mais on n'a pas \'egalit\'e dans le cas g\'en\'eral. On a cependant l'\'egalit\'e suivante comme dans le lemme \ref{mathfra}:
$$\hat{\Pi}(\hat{\mathfrak{K}_t})=\underset{i\in I}{\bigcap}\mathfrak{K}_{t,i}'$$
et on notera $\mathfrak{K}_t'$ cette sous-partie de $\Yim_t$.

\subsection{Les op\'erateurs $S_{t,i}^l$}

\subsubsection{D\'efinition}

On consid\`ere les $\Yim_t$-modules libres:
$$\Yim^l_{t,i}=\underset{a\in\CC^*}{\bigoplus}\Yim_t.S_{i,a}$$
On d\'eduit respectivement de $\hat{\Pi}_t$, $\hat{\Pi}$ (dans le cas $ADE$), $\Pi_t$ des applications $\hat{\Pi}_{t,i}^l$, $\hat{\Pi}_{i}^l$, $\Pi_{t,i}^l$.

\noindent On note $S_{t,i}^l$ l'application $\ZZ[t,t^{-1}]$-lin\'eaire $S_{t,i}^l:\Yim_t\rightarrow \Yim_{t,i}^l$ qui prend sur un mon\^ome $m\in\Yim_t$ la valeur:
$$S_{t,i}^l(m)=m(\underset{a\in\CC^*/u_{i,a}(m)\geq 0}{\sum}(1+...+t^{2(u_{i,a}(m)-1)})S_{i,a}-\underset{a\in\CC^*/u_{i,a}(m)<0}{\sum}(t^{-2}+...+t^{2u_{i,a}(m)})S_{i,a})$$
On voit imm\'ediatement:
$$S_{t,i}^l(Y_{j,a})=\delta_{j,i}Y_{i,a}S_{i,a}\text{ et }S_{t,i}^l(Y_{j,a}^{-1})=-\delta_{j,i}t^{-2}Y_{i,a}^{-1}S_{i,a}$$

\begin{lem}
Le diagramme (1) suivant est commutatif:
$$\begin{array}{rcccl}
\hat{\Yim}_t&\stackrel{\hat{S}_{t,i}^l}{\longrightarrow}&\hat{\Yim}_{t,i}^l\\
\hat{\Pi}_t\downarrow &&\downarrow&\hat{\Pi}_{t,i}^l\\
\Yim_t&\stackrel{S_{t,i}^l}{\longrightarrow}&\Yim_{t,i}^l\\
\Pi_t\downarrow &&\downarrow&\Pi_{t,i}^l\\
\Yim&\stackrel{S_i^l}{\longrightarrow}&\Yim_{1,i}^l\\
\end{array}$$ 
Dans le cas $ADE$, le diagramme (1)' obtenu en utilisant respectivement $\hat{\Pi},\hat{\Pi}_{i}^l$ \`a la place de $\hat{\Pi}_t,\hat{\Pi}_{t,i}^l$ est commutatif \'egalement.
\end{lem}

\demo

 Toutes les applications sont $\ZZ$-lin\'eaires, il suffit donc de regarder un mon\^ome $m\in\hat{\Yim}_t$ et $\lambda\in\ZZ[t,t^{-1}]$:
$$\begin{array}{ll}&(\hat{\Pi}_{t,i}^l\circ \hat{S}^l_{t,i})(\lambda m)
\\=&\lambda\hat{\Pi}_{t,i}^l(m(\underset{a\in\CC^*/u_{i,a}\geq 0}{\sum}(1+...+t^{2(u_{i,a}-1)})S_{i,a}-\underset{a\in\CC^*/u_{i,a}<0}{\sum}(t^{-2}+...+t^{2u_{i,a}})S_{i,a})
\\=&\lambda\hat{\Pi}_t(m)(\underset{a\in\CC^*/u_{i,a}\geq 0}{\sum}(1+...+t^{2(u_{i,a}-1)})S_{i,a}-\underset{a\in\CC^*/u_{i,a}<0}{\sum}(t^{-2}+...+t^{2u_{i,a}})S_{i,a})
\\=&\lambda\underset{a\in\CC^*}{\prod}Y_{i,a}^{u_{i,a}}(\underset{a\in\CC^*/u_{i,a}\geq 0}{\sum}(1+...+t^{2(u_{i,a}-1)})S_{i,a}-\underset{a\in\CC^*/u_{i,a}<0}{\sum}(t^{-2}+...+t^{2u_{i,a}})S_{i,a})
\\=&S_{t,i}^l(\lambda\underset{a\in\CC^*}{\prod}Y_{i,a}^{u_{i,a}})
\\=&S_{t,i}^l(\hat{\Pi}_t(\lambda m))\end{array}$$

\noindent puis un mon\^ome $m\in\Yim_t$ et $\lambda(t)\in\ZZ[t,t^{-1}]$:
$$\begin{array}{ll}&(\Pi_{t,i}^l\circ S^l_{t,i})(\lambda(t) m)
\\=&\lambda(1)\Pi_{t,i}^l(m(\underset{a\in\CC^*/u_{i,a}\geq 0}{\sum}(1+...+t^{2(u_{i,a}-1)})S_{i,a}-\underset{a\in\CC^*/u_{i,a}<0}{\sum}(t^{-2}+...+t^{2u_{i,a}})S_{i,a})
\\=&\lambda(1)m(\underset{a\in\CC^*}{\sum}u_{i,a}S_{i,a})
\\=&S_i^l(\lambda(1)m)
\\=&S_i^l(\hat{\Pi}_t(\lambda(t) m))\end{array}$$
Le diagramme (1)' se traite de mani\`ere analogue.
\qed

\subsubsection{Interpr\'etation des $S^l_i$ en terme de d\'erivation}

\begin{lem} Il existe sur $\Yim^l_{t,i}$ une unique structure de bimodule pour la multiplication usuelle $.$ de $\Yim_t$ telle la structure \`a gauche soit la structure ci-dessus, et que pour tout $a\in\CC^*$ et tout mon\^ome $\underset{a\in\CC^*,j\in I}{\prod}Y_{j,a}^{u_{j,a}}=m\in\Yim_t$:
$$S_{i,a}.m=t^{2u_{i,a}(m)}m.S_{i,a}\text{ , }S_{i,a}.t=t.S_{i,a}$$
\end{lem}

\noindent La d\'emonstration est compl\`etement analogue au cas $\hat{\Yim}^l_{t,i}$.

\noindent Notons qu'on a alors pour tout $a\in\CC^*$:
$$S_{i,a}.Y_{i,a}=t^2Y_{i,a}.S_{i,a}\text{ , }S_{i,a}.Y_{i,a}^{-1}=t^{-2}Y_{i,a}^{-1}.S_{i,a}$$

\begin{prop}
L'application $S_{t,i}^l$ a une propri\'et\'e de d\'erivation:
$$\forall U,V\in \hat{\Yim}_t,\hat{S}_{t,i}^l(U.V)=U.\hat{S}_{t,i}^l(V)+\hat{S}_{t,i}^l(U).V$$
C'est de plus l'unique d\'erivation $\ZZ[t,t^{-1}]$-lin\'eaire telle que $S_{t,i}^l(Y_a)=Y_a.S_a$.
\end{prop}

La d\'emonstration est compl\`etement analogue au cas $\hat{S}^l_{t,i}$.

\subsection{t-analogues des op\'erateurs d'\'ecrantage pour $\Yim_t$}

\subsubsection{D\'efinition des op\'erateurs $S_{t,i}$}

On consid\`ere le sous $\ZZ[t,t^{-1}]$-module $\hat{\Pi}_{t,i}^l(\hat{F}_{t,i})$ de $\Yim_{t,i}^l$. 

\begin{lem} Les \'el\'ements de $\hat{\Pi}_{t,i}^l(\hat{F}_{t,i})$ sont les \'el\'ements de $\Yim_{t,i}^l$ de la forme:
$$\underset{a\in\CC^*,m\in \hat{\Pi}_t(\hat{A})}{\sum}\lambda_{m,a}m(A_{i,aq_i}^{-1}t^{2u_{i,aq_i^2}(m)}S_{i,aq_i^2}-t^2S_{i,a})$$
avec les $\lambda_{m,a}\in \ZZ[t,t^{-1}]$ presque tous nuls. 
\end{lem}

\demo

Le r\'esultat d\'ecoule du fait que $\hat{\Pi}_t$ est un morphisme d'anneaux qui conserve les quantit\'es $u_{i,a}$.\qed

\noindent Soit \`a pr\'esent:
$$F_{t,i}=\underset{a\in\CC^*,m\in A}{\sum}\ZZ[t,t^{-1}].m(A_{i,aq_i}^{-1}t^{2u_{i,aq_i^2}(m)}S_{i,aq_i^2}-t^2S_{i,a})$$
C'est un sous $\ZZ[t,t^{-1}]$-module de $\Yim_{t,i}^l$ qui contient $\hat{\Pi}_{t,i}^l(\hat{F}_{t,i})$.
 
\noindent Notons que les \'el\'ements de $F_{t,i}$ sont les \'el\'ements de $\Yim_{t,i}^l$ de la forme:
$$\underset{a\in\CC^*}{\sum}(A_{i,aq_i}^{-1}S_{i,aq_i^2}.U_a-t^2U_a.S_{i,a})$$
avec les $U_a\in\Yim_t$.

\begin{defi} On appelle $\Yim_{t,i}$ le $\ZZ[t,t^{-1}]$-module quotient de $\Yim_{t,i}^l$ par $F_{t,i}$, et $S_{t,i}$ l'application obtenue \`a partir de $S^l_{t,i}$ par projection sur $\Yim_{t,i}$.\end{defi}

\noindent Notons que $F_{t,i}$ n'est pas un $\Yim_t$-sous module de $\Yim_{t,i}^l$, mais c'est une partie de $\Yim_{t,i}^l$ stable par multiplication par des \'el\'ements de $\ZZ[Y_{j,a}^{\pm}]_{j\neq i}$. En particulier on peut d\'efinir une multiplication \`a gauche sur $\Yim_{t,i}$ par les \'elements de $\ZZ[Y_{j,a}^{\pm}]_{j\neq i}$, qui commute avec la projection $p_i$ de $\Yim_{t,i}^l$ sur $\Yim_{t,i}$.

\begin{lem}
\noindent Les applications $\Pi_{t,i}^l,\hat{\Pi}_{t,i}^l$ donnent naturellement des applications $\Pi_{t,i},\hat{\Pi}_{t,i}$ rendant le diagramme (2) suivant commutatif:
$$\begin{array}{rcccl}
\hat{\Yim}_{t,i}^l&\longrightarrow&\hat{\Yim}_{1,i,t}\\
\hat{\Pi}_{t,i}^l\downarrow &&\downarrow&\hat{\Pi}_{t,i}\\
\Yim_{t,i}^l&\longrightarrow&\Yim_{t,i}\\
\Pi_{t,i}^l\downarrow &&\downarrow&\Pi_{t,i}\\
\Yim_{i}^l&\longrightarrow&\Yim_{i}\\
\end{array}$$
\end{lem}

\demo

Il suffit de v\'erifier que les applications $\ZZ$-lin\'eaires $\Pi_{t,i}^l,\hat{\Pi}_{t,i}^l$ passent aux quotient. Or $\hat{\Pi}_{t,i}^l(\hat{F}_{t,i})\subset F_{t,i}$, donc $\hat{\Pi}_{t,i}$ est bien d\'efinie. Puis pour $x\in F_{t,i}$, on a avec les notations d\'ej\`a utilis\'ees:
$$\Pi_{t,i}(x)=\underset{a\in\CC^*}{\sum}U_a(1)(A_{i,aq_i}^{-1}.S_{i,aq_i^2}-S_{i,a})\in F_i$$
\qed

\begin{prop} On a le diagramme commutatif suivant:
$$\begin{array}{rcccl}
\hat{\Yim}_t&\stackrel{\hat{S}_{t,i}}{\longrightarrow}&\hat{\Yim}_{t,i}\\
\hat{\Pi}_t\downarrow &&\downarrow&\hat{\Pi}_{t,i}\\
\Yim_t&\stackrel{S_{t,i}}{\longrightarrow}&\Yim_{t,i}\\
\Pi_t\downarrow &&\downarrow&\Pi_{t,i}\\
\Yim&\stackrel{S_i}{\longrightarrow}&\Yim_{i}\\
\end{array}$$ 
et on a :
$$\mathfrak{K}_{t,i}\subset \text{Ker}(S_{t,i})$$
\end{prop}

\demo 

La commutativit\'e du diagramme provient de la commutativit\'e des diagrammes (1) et (2). 

\noindent Puis $\hat{\Pi}_{t,i}\circ \hat{S}_{t,i}=S_{t,i}\circ \hat{\Pi}_t$ et $\hat{\mathfrak{K}}_{t,i}\subset \text{Ker}(\hat{S}_{t,i})$ implique $\hat{\Pi}_t(\hat{\mathfrak{K}}_{t,i})\subset \text{Ker}(S_{t,i})$. 

\noindent Soit alors $m\in B_i$ qu'on d\'ecompose $m=m_i\underset{j\neq i}{\prod}m_j$ avec les $m_j\in \ZZ[Y_{j,a}^{\pm}]_{a\in\CC^*}$. Alors:
$$E_{0,i}(m)=E_{0,i}(m_i)\underset{j\neq i}{\prod}m_j$$
et comme pour tout $a\in\CC^*$, $u_{i,a}(\underset{j\neq i}{\prod}m_j)=0$, on a:
$$S_{t,i}^l(E_{0,i}(m))=(\underset{j\neq i}{\prod}m_j)S_{t,i}^l(E_{0,i}(m_i))$$
Alors pour la multiplication \`a gauche sur $\Yim_{t,i}$ par des \'el\'ements de $\ZZ[Y_{j,a}^{\pm}]_{j\neq i}$, on a:
$$S_{t,i}(E_{0,i}(m))=(\underset{j\neq i}{\prod}m_j)S_{t,i}(E_{0,i}(m_i))$$
Mais alors comme $m$ est $i$-dominant, on a:
$$m_i=\underset{a\in\CC^*}{\prod}Y_{i,a}^{u_{i,a}}=\hat{\Pi}_t(\underset{a\in\CC^*}{\prod}W_{i,a}^{u_{i,a}})$$
avec les $u_{i,a}=u_{i,a}(m_i)\geq 0$. En cons\'equence: 
$$E_{0,i}(m_i)=\hat{\Pi}_t(E_i(\underset{a\in\CC^*}{\prod}W_{i,a}^{u_{i,a}}))\in \hat{\Pi}_t(\hat{\mathfrak{K}}_{t,i})\subset \text{Ker}(S_{t,i})$$\qed

\subsubsection{Remarques sur les op\'erateurs $S_{t,i}'$}

On peut faire une contruction analogue relative \`a $\mathfrak{K}'$ en consid\'erant les:
$$F_{t,i}'=\underset{a\in\CC^*,m\in A}{\sum}\ZZ[t,t^{-1}].m(A_{i,aq_i}^{-1}t^{u_{i,aq_i^2}(m)-u_{i,a}(m)}S_{i,aq_i^2}-tS_{i,a})$$
puis $\Yim_{t,i}'=\Yim_{t,i}^l/F_{t,i}'$, et $S_{t,i}'$ la compos\'ee de $S_{t,i}^l$ avec la projection de $\Yim_{t,i}^l$ sur $\Yim_{t,i}'$.

\noindent Dans le cas $ADE$, on a $\hat{\Pi}_{t,i}^l(\hat{F}_{t,i})\subset F_{t,i}'$, le diagramme commutatif:
$$\begin{array}{rcccl}
\hat{\Yim}_t&\stackrel{\hat{S}_{t,i}}{\longrightarrow}&\hat{\Yim}_{t,i}\\
\hat{\Pi}\downarrow &&\downarrow&\hat{\Pi}_{i}\\
\Yim_t&\stackrel{S_{t,i}'}{\longrightarrow}&\Yim_{t,i}'\\
\Pi_t\downarrow &&\downarrow&\Pi_{t,i}'\\
\Yim&\stackrel{S_i}{\longrightarrow}&\Yim_{i}\\
\end{array}$$ 
et on a: 
$$\mathfrak{K}_{t,i}'\subset \text{Ker}(S_{t,i}')$$

\subsection{Noyau des $t$-op\'erateurs d'\'ecrantage $S_{t,i}$}

\begin{thm}\label{ksimple} On a $\mathfrak{K}_{t,i}= \text{Ker}(S_{t,i})$.\end{thm}

\noindent On pourrait montrer ce r\'esultat de la m\^eme mani\`ere que $\hat{\mathfrak{K}}_{t,i}= \text{Ker}(\hat{S}_{t,i})$ en utilisant la d\'ecomposition de $\Yim_t$ du lemme \ref{sdirecte}. Cette m\'ethode permet aussi retrouver le r\'esultat du th\'eor\`eme \ref{frene} en utilisant la d\'ecomposition de $\Yim$:
$$\Yim=\mathfrak{K}_i\oplus \underset{m\in A-B_i}{\bigoplus} \ZZ m$$
On propose ici une alternative qui d\'eduit le r\'esultat du th\'eor\`eme \ref{ksimple} de celui du th\'eor\`eme \ref{frene}. Elle n\'ecessite quelques lemmes pr\'eliminaires.

\noindent Noter que tout ce qui suit peut \^etre appliqu\'e de mani\`ere analogue \`a $S_{t,i}'$ dans le cas $ADE$, ce qui donne $\mathfrak{K}_{t,i}'= \text{Ker}(S_{t,i}')$.

\subsubsection{Lemmes pr\'eliminaires}

\begin{lem}\label{decomp} Tout $u\in \Yim_t$ s'\'ecrit de mani\`ere unique sous la forme:
$$u=\underset{m\in D}{\sum}(t-1)^{p(m)}t^{-q(m)}(\alpha_0(m)+\alpha_1(m)(t-1)+\alpha_2(m)(t-1)^2+...)m$$
avec $D\subset A$, les $p(m),q(m),\alpha_0(m),\alpha_1(m),...\in\NN$ et $\alpha_0(m)\neq 0$.
\end{lem}

\demo 

On d\'ecompose $u$ sur la somme directe $\Yim_t=\underset{m\in A}{\sum} \ZZ[t,t^{-1}]m$, et il suffit donc de consid\'erer un polyn\^ome de Laurent $P\in\ZZ[t,t^{-1}]$ non nul et de montrer qu'il s'\'ecrit de mani\`ere unique: 
$$P=\frac{(t-1)^{p(m)}}{t^{q(m)}}(\alpha_0(m)+\alpha_1(m)(t-1)+\alpha_2(m)(t-1)^2+...)$$ 
Si $P\in\ZZ[t]$, c'est le cas car on a une base gradu\'ee $((t-1)^p)_p$ de $\ZZ[t]$. Dans le cas g\'en\'eral, $P$ s'\'ecrit de mani\`ere unique $P=t^{-q(m)}Q$ avec $Q\in\ZZ[t]$ et $q(m)\in\NN$.\qed

\begin{cor}\label{noypit} Le noyau de $\Pi_t$ est $\text{Ker}(\Pi_t)=(t-1)\Yim_t$.\end{cor}

\demo

L'inclusion $(t-1)\Yim_t\subset \text{Ker}(\Pi_t)$ est claire, puis si $u\in \text{Ker}(\Pi_t)$, en utilisant la d\'ecomposition du lemme \ref{decomp}, on voit que: 
$$\underset{m\in D,p(m)=0}{\sum}\alpha_0(m)m=0$$
or comme les $\alpha_0(m)$ sont non nuls, pour $m\in D$ on a $p(m)>0$. \qed

\begin{lem}\label{liberte} Si $\alpha\in \Yim_{t,i}^l$ v\'erifie $(t-1)\alpha\in F_{t,i}$ alors $\alpha\in F_{t,i}$.\end{lem}

\demo

Pour un tel  $\alpha\in \Yim_{t,i}^l$, on peut \'ecrire:
$$(t-1)\alpha=\underset{a\in\CC^*,m\in A}{\sum}\lambda_{m,a}m(A_{i,aq_i}^{-1}t^{2u_{i,aq_i^2}(m)}S_{i,aq_i^2}-t^2S_{i,a})$$
avec les $\lambda_{m,a}\in \ZZ[t,t^{-1}]$ presque tous nuls. Mais si on \'evalue cette expression \`a $t=1$, on trouve dans $\Yim^l_{1,i}$:
$$0=\underset{a\in\CC^*,m\in A}{\sum}\lambda_m(1)m(A_{i,aq_i}^{-1}S_{i,aq_i^2}-S_{i,a})=\underset{a\in\CC^*}{\sum}U_a(1)(A_{i,aq_i}^{-1}S_{i,aq_i^2}-S_{i,a})$$
avec:
$$U_a=\underset{m\in A}{\sum}\lambda_{m,a}m$$ 
Supposons alors par l'absurde qu'il existe un $a$ tel que $U_a(1)\neq 0$. On consid\`ere alors la plus grande puissance de $q$ tel que $U_{aq^m}(1)\neq 0$ (qui existe car les $U_b$ sont presque tous nuls). Alors $A_{i,aq^{m+1}}^{-1}U_{aq^m}(1)$ est le coefficient de $S_{i,aq_i^2}$, donc $A_{i,aq^{m+1}}^{-1}U_{aq^m}(1)=0$, contradiction. On peut donc \'ecrire tous les $U_a$ sous la forme $U_a=(t-1)U_a'$ avec les $U_a'=\underset{m\in A}{\sum}\lambda_{m,a}'m\in \Yim_t$, et $(t-1)\alpha=(t-1)\beta$ avec:
$$\beta=\underset{a\in\CC^*,m\in A}{\sum}\lambda_{m,a}m(t^{2u_{i,aq_i^2}(m)}A_{i,aq_i}^{-1}.S_{i,aq_i^2}-t^2S_{i,a})\in F_{t,i}$$
Mais comme $\Yim_{t,1}^l$ est un $\ZZ[t,t^{-1}]$-module libre, on a $\alpha=\beta$. \qed

\subsubsection{D\'emonstration du th\'eor\`eme \ref{ksimple}}

\demo

La premi\`ere inclusion $\mathfrak{K}_{t,i}\subset \text{Ker}(S_{t,i})$ est d\'ej\`a connue.

\noindent Consid\'erons $u\in \text{Ker}(S_{t,i})$. Alors $0=\Pi_{t,i}(S_{t,i}(u))=S_i(\Pi_t(u))$ et donc d'apr\`es le th\'eor\^eme \ref{frene}, $\Pi_t(u)\in\ZZ[Y_{j,b}^{\pm}]_{j\neq i,b\in\CC^*}\ZZ[Y_{i,a}(1+A_{i,aq_i}^{-1})]_{a\in\CC^*}$, soit:
$$\Pi_t(u)=\underset{m}{\sum}\lambda_m\underset{a\in\CC^*}{\prod}(Y_{i,a}(1+A_{i,aq_i}^{-1}))^{p_a(m)}$$
avec $\lambda_m\in \ZZ[Y_{j,b}^{\pm}]_{j\neq i,b\in\CC^*}$ et les $p_a(m)\geq 0$. Pour chaque $m$, on pose:
$$v_m=\underset{a\in\CC^*}{\prod}(W_{i,a}(1+V_{i,aq_i}^{-1}))^{*p_a(m)}=E_i(\underset{a\in\CC^*}{\prod}W_{i,a}^{p_a(m)})\in \hat{\mathfrak{K}}_{t,i}$$
En consid\'erant alors:
$$v=\underset{m}{\sum}\lambda_m\hat{\Pi}_t(v_m)$$
on a $v\in \mathfrak{K}_{t,i}$ et $\Pi_t(v)=\Pi_t(u)$, donc $v-u\in \text{Ker}(\Pi_t)=(t-1)\Yim_t$ d'apr\`es le corollaire \ref{noypit}. En cons\'equence:
$$u=v+(t-1)u_1$$
avec $u_1\in\Yim_t$. Mais alors $(t-1)u_1\in \text{Ker}(S_{t,i})$, soit $(t-1)S_{t,i}^l(u_1)\in F_{t,i}$. Alors d'apr\`es le lemme \ref{liberte}, $S_{t,i}^l(u_1)\in F_{t,i}$, soit $u_1\in \text{Ker}(S_{t,i})$. On peut recommencer avec $u_1$, et on obtient par r\'ecurrence que pour tout $p\in\NN$, il existe $w_p\in \mathfrak{K}_{t,i}$ et $u_p\in \text{Ker}(S_{t,i})$ tels que:
$$u=w_p+(t-1)^pz_p$$ 
D\'ecomposons $u=b+c$ sur la somme directe du lemme \ref{sdirectesimp}:
$$\Yim_t=\mathfrak{K}_{t,i}\oplus \underset{m\in \hat{A}-\hat{B}_i}{\bigoplus} \ZZ[t,t^{-1}]m$$
 et supposons par l'absurde que $c\neq 0$. Pour $p$, prenons $p_0=p(m_0)+1$ avec $p(m_0)$ le plus grand $p(m)$ qui apparait dans la d\'ecomposition de $c$ du lemme \ref{decomp}. On obtient une \'ecriture $u=w+(t-1)^{p_0}v$.
D\'ecomposons $v=b'+c'$ sur la somme directe du lemme \ref{sdirecte}. Alors:
$$u=b+c=w+(t-1)^{p_0}b'+(t-1)^{p_0}c'$$ 
et $b=w+(t-1)^{p_0}b'$ et $c=(t-1)^{p_0}c'$. Donc les $p(m)$ qui apparaissent dans la d\'ecomposition de $c$ du lemme \ref{decomp} sont tous strictement plus grands que $p_0$, contradiction. On a donc $\mathfrak{K}_{t,i}=\text{Ker}(S_{t,i})$.
\qed

\section{Compl\'ements relatifs aux involutions}

On rappelle les involutions d\'efinies par Nakajima: sur $\Yim_t$ on pose $\overline{t}=t^{-1},\overline{Y_{i,a}^{\pm}}=Y_{i,a}^{\pm}$, et sur $\hat{\Yim}_t$ pour $d$ un bicarat\`ere, on pose:
$$\overline{t}=t^{-1}\text{ , }\overline{m}=t^{2d(m,m)}m$$
En particulier dans le cas $ADE$ on a l'involution $\overline{ }$ obtenue avec $d_N$. Elle est alors anti multiplicative relativement \`a $*$ et commute avec $\hat{\Pi}$.

\noindent On \'etend ces involutions \`a $\Yim_{t,i}^l$ (respectivement \`a $\hat{\Yim}_{t,i}^l$) en posant $\overline{S_{i,a}}=t^{-2}S_{i,a}$, soit:
$$\overline{\underset{a\in\CC^*}{\sum}U_aS_{i,a}}=\underset{a\in\CC^*}{\sum}t^{-2}S_{i,a}\overline{U_a}$$
pour des $U_a$ dans $\Yim_t$ (respectivement $\hat{\Yim}_t$).

\begin{lem} On a pour $x\in\hat{\Yim}_t, y\in\Yim_t$:
$$\overline{S_{t,i}^l(x)}=S_{t,i}^l(\overline{x})\text{ , }\overline{\hat{S}_{t,i}^l(y)}=\hat{S}_{t,i}^l(\overline{y})$$
De plus $\hat{F}_{t,i},F_{t,i}'$ sont stables par les involutions correpondantes.
\end{lem}

\demo

Les deux r\'esultats s'obtiennent de mani\`ere analogue, en consid\'erant par exemple un mon\^ome $a\in\Yim_t$:
$$\begin{array}{ll}&\overline{S_{t,i}^l(\lambda(t)m)}
\\=&\lambda(t^{-1})(\underset{a\in\CC^*u_{i,a}(m)\geq 0}{\sum}(1+t^{-2}+...+t^{-2(u_{i,a}(m)-1)})t^{-2}S_{i,a}m
\\&+\underset{a\in\CC^*u_{i,a}(m)<0}{\sum}(t^{2}+...+t^{-2u_{i,a}(m)})t^{-2}S_{i,a}m)
\\=&\lambda(t^{-1})m(\underset{a\in\CC^*u_{i,a}(m)\geq 0}{\sum}(t^{2(u_{i,a}(m)-1)}+...+t^2+1)S_{i,a}+\underset{a\in\CC^*u_{i,a}(m)<0}{\sum}(t^{2u_{i,a}(m)}+...+t^{-2})S_{i,a})
\\=&\lambda(t^{-1})S_{t,i}^l(m)=S_{t,i}^l(\overline{\lambda(t)m})\end{array}$$
Consid\'erons ensuite, par exemple dans le cas $ADE$:
$$\begin{array}{ll}&\overline{\lambda(t)m(V_{i,aq_i}t^{2u_{i,aq_i^2}(m)}S_{i,aq_i^2}-t^2S_{i,a})}
\\=&\lambda(t^{-1})t^{-2}(t^{-2u_{i,aq_i^2}(m)}S_{i,aq_i^2}t^{2d(m,m)-2+2u_{i,a}(m)+2u_{i,aq_i^2}(m)}mV_{i,aq_i}-t^{-2}S_{i,a}t^{2d(m,m)}m)
\\=&\lambda(t^{-1})t^{-2+2d(m,m)}m(V_{i,aq_i}t^{-2}S_{i,aq_i^2}t^{-2+2u_{i,a}(m)+2u_{i,aq_i^2}(m)}-t^{-2+2u_{i,a}(m)}S_{i,a})
\\=&\lambda(t^{-1})t^{-6+2d(m,m)+2u_{i,a}(m)}m(V_{i,aq_i}S_{i,aq_i^2}t^{2u_{i,aq_i^2}(m)}-t^2S_{i,a})\in\hat{F}_{t,i}\end{array}$$
et dans le cas g\'en\'eral:
$$\begin{array}{ll}&\overline{\lambda(t)m(A_{i,aq_i}^{-1}t^{u_{i,aq_i^2}(m)-u_{i,a}(m)}S_{i,aq_i^2}-tS_{i,a})}
\\=&\lambda(t^{-1})t^{-2}(t^{u_{i,a}(m)-u_{i,aq_i^2}(m)}S_{i,aq_i^2}mA_{i,aq_i}^{-1}-t^{-1}S_{i,a}m)
\\=&\lambda(t^{-1})t^{-2}m(A_{i,aq_i}^{-1}t^{-2+u_{i,aq_i^2}(m)+u_{i,a}(m)}S_{i,aq_i^2}-t^{-1+2u_{i,a}(m)}S_{i,a})
\\=&\lambda(t^{-1})t^{-4+2u_{i,a}(m)}m(A_{i,aq_i}^{-1}t^{u_{i,aq_i^2}(m)-u_{i,a}(m)}S_{i,aq_i^2}-tS_{i,a})\in F_{t,i}'\end{array}$$
\qed

\noindent On peut ainsi d\'efinir des involutions sur $\hat{\Yim}_{t,i},\Yim_{t,i}'$ qui commutent respectivement avec les op\'erateurs d'\'ecrantage associ\'es.

\noindent {\bf Remerciements}: Je remercie M. Rosso pour nos discussions et ses pr\'ecieux conseils, et H. Nakajima pour ses indications sur les $q,t$-caract\`eres.


\begin{thebibliography}{99}

\bibitem[1]{bou} {\bf N. Bourbaki}, {\it Groupes et alg\`ebres de Lie} 

{Chapitres IV-VI, Hermann (1968)}
\mk

\bibitem[2]{Fre} {\bf E. Frenkel et N. Reshetikhin}, {\it The $q$-Characters of Representations of Quantum Affine Algebras and Deformations of $W$-Algebras} 

{http://www.arxiv.org/abs/math/9810055}

{Recent Developments in Quantum Affine Algebras and related topics, Cont. Math., vol 248, pp 163-205 (1999)}
\mk

\bibitem[3]{Fre2} {\bf E. Frenkel et E. Mukhin}, 
{\it Combinatorics of $q$-Characters of Finite-Dimensional Representations of Quantum Affine Algebras} 

{http://www.arxiv.org/abs/math/9911112}

{Comm. in Math. Phy., vol 216, no. 1, pp 23-57 (2001)}

\mk

\bibitem[4]{Naa} {\bf H. Nakajima}, {\it $t$-Analogue of the $q$-Characters of Finite Dimensional Representations of Quantum Affine Algebras} 

{http://www.arxiv.org/abs/math/0009231} 

``Physics and Combinatorics'', Proc. Nagoya 2000 International Workshop, World Scientific, pp 181-212 (2001)

\mk

\bibitem[5]{Nab} {\bf H. Nakajima}, {\it Quiver Varieties and $t$-Analogs of $q$-Characters of Quantum Affine Algebras} 

{http://www.arxiv.org/abs/math/0105173} 

\mk

\bibitem[6]{Ro} {\bf M. Rosso}, 
{\it Repr\'esentations des groupes quantiques} 

{S\'eminaire Bourbaki exp. no. 744, Ast\'erisque 201-203, 443-83, SMF (1992)}

\end{thebibliography}
\end{document}